\documentclass[11pt]{article}
\usepackage[dvips]{graphics}
\usepackage{multirow}
\usepackage{array}
\usepackage{epsfig,rotating}
\usepackage{amssymb,amsmath}
\usepackage{latexsym}
\usepackage[usenames]{color}

\numberwithin{equation}{section}

\newcommand{\be}{\begin{equation}}
\newcommand{\ee}{\end{equation}}
\newcommand{\beaa}{\begin{eqnarray*}}
\newcommand{\eeaa}{\end{eqnarray*}}
\newcommand{\bea}{\begin{eqnarray}}
\newcommand{\eea}{\end{eqnarray}}
\newcommand{\lbl}{\label}

\newcommand{\bei}{\begin{itemize}}
\newcommand{\eei}{\end{itemize}}

\newcommand{\bd}{\bold}

\newtheorem{theorem}{ \noindent T{\footnotesize HEOREM}}
\newtheorem{prop}{ \noindent P{\footnotesize ROPOSITION}}[section]
\newtheorem{lemma}{ \noindent L{\footnotesize EMMA}}[section]
\newtheorem{coro}{ \noindent C{\footnotesize OROLLARY}}
\begin{document}

%Coherence-Random-Matrix02092011
%\title{On the Coherence of High-Dimensional Random Matrices}
\title{Central Limit Theorems for Classical Likelihood Ratio Tests for High-Dimensional Normal Distributions}
\author{Tiefeng Jiang$^1$ and Fan Yang$^{1,2}$\\
University of Minnesota}

\date{}
\maketitle

\footnotetext[1]{School of Statistics, University of Minnesota, 224 Church
Street, S.E., MN55455, USA, jiang040@umn.edu. \newline  \indent \  \ The research of Tiefeng Jiang was
supported in part by NSF FRG Grant DMS-0449365 and NSF Grant DMS-1208982.}
\footnotetext[2]{Boston Scientific, 1 Scimed Place, Maple Grove, MN 55311, USA,
%\newline  \indent \
yangf@bsci.com.}

\begin{abstract}
\noindent For random samples of size $n$ obtained from $p$-variate normal distributions, we
consider the classical likelihood ratio tests (LRT) for their means and covariance matrices
in the high-dimensional setting.
These test statistics have been extensively studied in
multivariate analysis and their limiting distributions under
the null hypothesis were proved to be chi-square distributions as $n$ goes to infinity and $p$ remains fixed. In this paper, we
consider the high-dimensional case where both $p$ and $n$ go to infinity with $p/n \to y \in (0, 1].$ We prove that the likelihood ratio test statistics under this
assumption will converge in distribution to normal distributions with explicit means and variances.
We perform the simulation study
to show that the likelihood ratio tests using our central limit theorems
outperform those using the traditional chi-square approximations for analyzing
high-dimensional data.

\end{abstract}

\noindent \textbf{Keywords:\/} Likelihood ratio test, central limit theorem, high-dimensional data, multivariate normal distribution, hypothesis test, covariance matrix, mean vector, multivariate Gamma function.

\noindent\textbf{AMS 2000 Subject Classification: \/} Primary 62H15; secondary 62H10.\\

%\textbf{Two things to consider:}

%1. we studied Gassian in this thesis, in practice our random variables are non-Gaussian, how sensitive for a small change?

%{\it Answer: our proof is based on mgf method, which has to rely on the Gaussian assumption. In the non-Gaussian case, the RMT theory help. For %example, Bai, Jiang , Yao and Zheng.}}

%\newpage

%\textbf{\large Contents}

%{\small
%\ref{Introduction}. Introduction

%%\ \ \ \ \ref{LLCR}. Limiting Laws of the Coherence of a Random Matrix

%\ \ \ \ \ref{Given_matrix}. Testing Covariance Matrices of Normal Distributions Proportional to Identity Matrix

%%\ \ \ \ \ref{iid.sec}. The i.i.d. Case

%%\ \ \ \ \ref{dependent.sec}. The Dependent Case

%\ \ \ \ \ref{Independence}. Testing Independence of Components of Normal Distributions

%\ \ \ \ \ref{Identical}. Testing that Several Normal Distributions Are Identical

%\ \ \ \ \ref{Constitution}. Testing Equality of Several Covariance Matrices

%\ \ \ \ \ref{Brown}. Testing Specified Values for Mean Vector and Covariance Matrix

%\ \ \ \ \ref{Correlation}. How Large Is the Determinant of A Sample Correlation Matrix?

%\ref{Proofs} Proofs

%\ \ \ \ \ref{preparation}. A Preparation

%\ \ \ \ \ref{Theorem_yet}. Proof of Theorem \ref{yet}

%\ \ \ \ \ref{Theorem_energy}. Proof of Theorem \ref{energy}

%\ \ \ \ \ref{Theorem_town}. Proof of Theorem \ref{town}

%\ \ \ \ \ref{Theorem_eve}. Proof of Theorem \ref{eve}

%\ \ \ \ \ref{Theorem_Ruth}. Proof of Theorem \ref{Ruth}

%\ \ \ \ \ref{Theorem_season}. Proof of Theorem \ref{season}

%%\ref{ladies} Appendix
%}

\newpage

\section{Introduction}\lbl{Introduction}

\noindent Traditional statistical theory, particularly in multivariate analysis, does not contemplate the demands of high dimensionality in data analysis due to  technological limitations and/or motivations. Consequently, tests of hypotheses and many other modeling procedures in many  classical textbooks of multivariate analysis such as Anderson (1958), Muirhead (1982), and Eaton (1983) are well developed under the assumption that the dimension of the dataset, denoted by $p$, is considered a fixed small constant or at least negligible compared with the sample size $n$. However, this assumption is no longer true for many modern datasets, because their dimensions can be proportionally large compared with the sample size. For example, the financial data,  the consumer data, the modern manufacturing data and the multimedia data all have this feature. More examples of high-dimensional data can be found in Donoho (2000) and Johnstone (2001).

Recently,   Bai et al. (2009) develop corrections to the traditional likelihood ratio test (LRT) to make it suitable for testing a high-dimensional normal distribution $N_p(\bd{\mu}, \bd{\Sigma})$ with $H_0:\  \bd{\Sigma} = \bd{I}_p\ \ \mbox{vs}\ \ H_a:\  \bd{\Sigma} \ne \bd{I}_p.$ The test statistic is chosen to be $L_n :=  \mbox{tr} (\bd{S}) -\log |\bd{S}| -p,$ where $\bd{S}$ is the sample covariance matrix from the data. In their derivation, the dimension $p$ is no longer considered a fixed constant, but rather a variable that goes to infinity along with the sample size $n$, and the ratio between $p=p_n$ and $n$ converges to a constant $y$, i.e.,
\begin{eqnarray}
\lim_{n \to \infty} \frac{p_n}{n} = y \in (0,1).
\end{eqnarray}
Jiang et al. (2012) further extend Bai's result to cover the case of $y=1.$
%and also proposed a new LRT for testing equality of two covariance matrices of normal distributions in high-dimensional situation.

In this paper, we study several other classical likelihood ratio tests for means and covariance matrices of high-dimensional normal distributions. Most of these tests have the asymptotic results for their test statistics derived decades ago under the assumption of a large $n$ but a fixed $p$. Our results supplement these traditional results in providing alternatives to analyze high-dimensional datasets including the critical case $p/n \to 1.$ We will briefly introduce these likelihood ratio tests next. In Section \ref{medical_summary}, for each LRT  described below, we first review the existing literature, then give our central limit theorem (CLT) results when the dimension and the sample size are comparable. We also make graphs and tables on the sizes and powers of these CLTs based on our simulation study to show that, as both $p$ and $n$ are large,  the traditional chi-square approximation behaves poorly and our CLTs improve the approximation very much.
\begin{itemize}

\item In Section \ref{Given_matrix}, for the normal distribution $N_p(\bd{\mu}, \bd{\Sigma}),$ we study the sphericity  test $H_0:\, \bd{\Sigma}=\lambda \bd{I}_p\ \ \mbox{vs}\ \ H_a:\, \bd{\Sigma} \ne \lambda \bd{I}_p$ with $\lambda$ unspecified. We derive the central limit theorem for the LRT statistic when $p/n \to y\in (0, 1].$ Its proof is given at Section \ref{Theorem_yet}.

\item In Section \ref{Independence}, we derive the CLT for the LRT statistic in testing that several components of a vector with distribution $N_p(\mu, \bd{\Sigma})$ are independent. The proof is presented at Section \ref{Theorem_energy}.

\item In Section \ref{Identical}, we consider the LRT  with $H_0: N_p(\mu_1, \bd{\Sigma}_1)=\cdots =N_p(\mu_k, \bd{\Sigma}_k),$ that is, several normal distributions are identical. We prove a CLT for the LRT statistic with the assumption $p/n_i\to y_i\in (0, 1]$ where $n_i$ is the sample size of a data set from $N_p(\mu_i, \bd{\Sigma}_i)$ for $i=1,2,\cdots,k.$ The proof of the theorem  is arranged at Section \ref{Theorem_town}.

\item In Section \ref{Constitution}, the test of the equality of the covariance matrices from several normal distributions are studied, that is, $H_0: \bd{\Sigma}_1=\cdots =\bd{\Sigma}_k.$ The LRT statistic is evaluated under the assumption $p/n_i\to y_i\in (0, 1]$ for $i=1,\cdots, k.$ This generalizes the work of Bai et al. (2009) and Jiang et al. (2012) from $k=2$ to any $k\geq 2.$ The proof of our result is given at Section \ref{Theorem_eve}.

\item In Section \ref{Brown}, we investigate LRT with $H_0: \bd{\mu}=\bd{0},\ \bd{\Sigma}=\bd{I}_p$ for the population distribution $N_p(\bd{\mu}, \bd{\Sigma}).$ With the dimension $p$ and the sample size $n$ satisfying $p/n\to y\in (0, 1],$ we derive the CLT for the LRT statistic. The corresponding theorem is proved at Section \ref{Theorem_Ruth}.

\item In Section \ref{Correlation}, we study the test that the population correlation matrix of a normal distribution is equal to an identity matrix, that is, all of the components of a normal vector  are independent (but not necessarily identically distributed). This is different from the test in Section \ref{Independence}  that {\it several components} of a normal vector are independent. The proof is presented at Section \ref{Theorem_season}.

\item In Sections \ref{simulation_study} and \ref{Conclusion_Discussion}, we show some simulation results,  state our method of the proofs and conclude by offering some open  problems.

\end{itemize}
One can see the value of $y=\lim (p/n)$ or $y_i=\lim (p/n_i)$ introduced above is restricted to  the range that $y\leq 1.$ In fact, when $y>1,$ some matrices involved in the LRT statistics do not have a full rank,  and consequently their determinants are equal to zero. As a result, the LRT statistics are not defined, or do not exist.

To our knowledge the central limit theorem of the LRT statistics mentioned above in the context of $p/n\to y\in (0, 1]$ are new in the literature. Similar research are Bai et al. (2009) and Jiang et al. (2012). The methods of the proofs in the three papers are different: the Random Matrix Theory is used in Bai et al. (2009); the Selberg integral is used in Jiang et al. (2012). Here we obtain the central limit theorems by analyzing the moments of the LRT statistics.

The organization of the rest of the paper is stated as follows. In Section \ref{medical_summary}, we give the details for each of the six tests described above. A simulation study on the sizes and powers of these tests is presented in Section \ref{simulation_study}. A discussion is given in Section \ref{Conclusion_Discussion}. The theorems appearing in each section are proved in Section \ref{Proofs}.  An auxiliary result on complex analysis is proved in the Appendix.

\section{Main Results}\lbl{medical_summary}
In this section we present the central limit theorems of six classical LRT statistics mentioned in the Introduction. The
six central limit theorems are stated in the following six subsections.

\subsection{Testing Covariance Matrices of Normal Distributions Proportional to Identity Matrix}\lbl{Given_matrix}

For distribution $N_p(\bd{\mu}, \bd{\Sigma}),$ we consider the spherical test
\bea\lbl{have}
H_0:\, \bd{\Sigma}=\lambda \bd{I}_p\ \ \mbox{vs}\ \ H_a:\, \bd{\Sigma} \ne \lambda \bd{I}_p
\eea
with $\lambda$ unspecified. Let $\bd{x}_1, \cdots, \bd{x}_n$ be i.i.d. $\mathbb{R}^p$-valued random variables with normal distribution $N_p(\bd{\mu}, \bd{\Sigma}).$  Recall
\begin{eqnarray}\lbl{golden}
\bar{\bd{x}}=\frac{1}{n}\sum_{i=1}^n\bd{x}_i\ \ \mbox{and}\ \ \bd{S}=\frac{1}{n}\sum_{i=1}^n(\bd{x}_i- \bar{\bd{x}})(\bd{x}_i- \bar{\bd{x}})'.
\end{eqnarray}
The likelihood ratio test statistic of  (\ref{have}) is first derived by Mauchly (1940) as
\begin{eqnarray}\lbl{Target}
V_n=|\bd{S}|\cdot \Big(\frac{\mbox{tr}(\bd{S})}{p}\Big)^{-p}.
\end{eqnarray}
By Theorem 3.1.2 and Corollary 3.2.19 from Muirhead (1982), under $H_0$ in (\ref{have}),
\bea\lbl{son}
\frac{n}{\lambda}\cdot \bd{S} \ \ \mbox{and}\ \  \bd{Z}'\bd{Z}\ \mbox{have the same distribution}
\eea
where $\bd{Z}:=(z_{ij})_{(n-1)\times p}$ and $z_{ij}$'s are i.i.d. with distribution $N(0, 1).$ This says that, with probability one,  $\bd{S}$ is not of full rank when $p \geq n$, and consequently $|\bd{S}|=0.$ This indicates that the likelihood
ratio test of (\ref{have}) only exists when $p \leq n-1$. The statistic $V_n$ is commonly known as
the ellipticity statistic. Gleser (1966) shows that the likelihood ratio test with the
rejection region $\{V_n \leq c_{\alpha}\}$ (where $c_{\alpha}$ is chosen so that the test has a significance
level of $\alpha$) is unbiased. A classical asymptotic  result shows that
\bea\lbl{afternoon}
-(n-1)\rho \log V_n\ \ \mbox{converges to}\ \ \chi^2(f)
\eea
in distribution as $n\to\infty$ with $p$ fixed, where
\bea\lbl{fee}
\rho=1-\frac{2p^2 + p +2}{6(n-1)p}\ \ \ \mbox{and}\ \ \ f=\frac{1}{2}(p-1)(p+2).
\eea
This can be seen from, for example, Theorem 8.3.7 from Muirhead (1982), the Slutsky lemma and the fact that $\rho=\rho_n\to 1$ as $n\to\infty$ and $p$ is fixed. The quantity $\rho$ is a correction term to improve the convergence rate.

Now we consider the case when both $n$ and $p$ are large. For clarity of taking limit, let $p=p_n$, that is, $p$ depends on $n.$

\begin{theorem}\lbl{yet} Let  $n>p+1$ for all $n\geq 3$ and $V_n$ be as in (\ref{Target}).   Assume $\lim_{n\to\infty}p/n= y\in (0, 1]$, then, under $H_0$ in (\ref{have}), $(\log V_n - \mu_n)/\sigma_n$ converges in distribution to  $N(0, 1)$ as $n\to\infty,$
%\begin{eqnarray*}
%\frac{\log V_n + p-\big(p-n+\frac{3}{2}\big)\log (1-\frac{p}{n-1})}{\sigma_n}\ \ \mbox{converges in %distribution to}\ \ N(0, 1)
%\end{eqnarray*}
%as $n\to\infty,$
where
\beaa
& & \mu_n=-p-\big(n-p-\frac{3}{2}\big)\log (1-\frac{p}{n-1})\ \ \mbox{and}\ \ \ \ \ \ \ \ \ \\
& & \sigma_n^2=-2\Big[\frac{p}{n-1}+\log (1-\frac{p}{n-1})\Big].
\eeaa
%where $\sigma_n^2=\frac{2p}{n-1}-2\log (1-\frac{p}{n-1}).$
\end{theorem}
As discussed below (\ref{son}), the LRT exists as $n\geq p+1$, however, we need a slightly stronger condition that $n>p+1$ because of the definition of $\sigma_n^2.$ Though $\lambda$ in (\ref{have}) is unspecified, the limiting distribution in Theorem \ref{yet} is pivotal, that is, it does not depend on $\lambda.$ This is because $\lambda$ is canceled in the expression of $V_n$ in (\ref{Target}):  $|\alpha \bd{S}|=\alpha^p|\bd{S}|$ and $(\mbox{tr}(\alpha \bd{S}))^{-p}=\alpha^{-p}\cdot (\mbox{tr}(\bd{S}))^{-p}$ for any $\alpha>0.$

\begin{figure}[t]
\centerline{\includegraphics[scale=0.44, bb=0 0 1000 525]{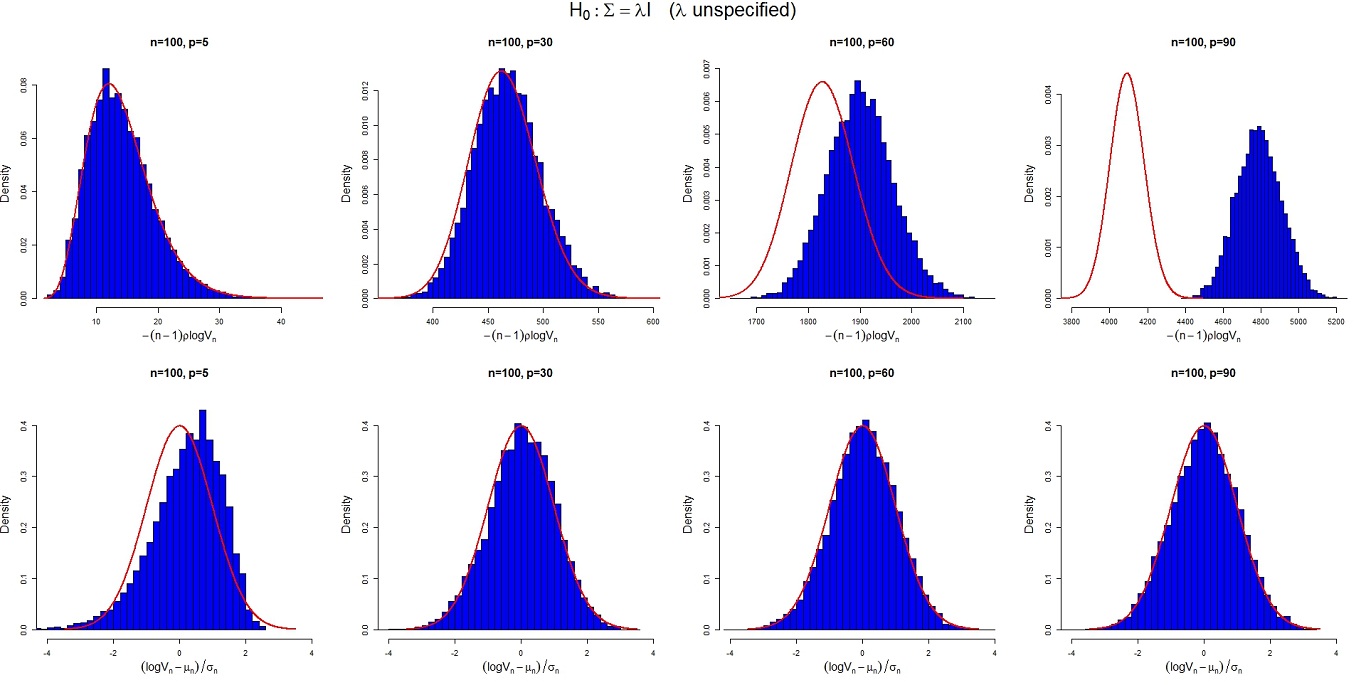}}
\caption{\sl Comparison between Theorem \ref{yet} and (\ref{afternoon}). We choose $n=100$ with $p=5, 30, 60, 90.$ The pictures in the top row show that the $\chi^2$ curves stay farther away from the histogram of $-(n-1)\rho\log V_n$  as $p$ grows. The bottom row shows that the $N(0, 1)$-curve fits the histogram of $(\log V_n-\mu_n)/\sigma_n$  better as $p$ becomes larger. }
%{\red the word ``$-(n-1)\rho \log V_n$" is too small. This is true for all pictures.}
\label{fig1}
\end{figure}

Simulation is run on the approximation in (\ref{afternoon}) and the CLT in Theorem \ref{yet}. The summary is given in Figure \ref{fig1}. It is seen from  Figure \ref{fig1} that the approximation in (\ref{afternoon}) becomes poorer as $p$ becomes larger relative to $n,$ and at the same time the CLT  in Theorem \ref{yet} becomes more precise. In fact, the chi-square approximation in (\ref{afternoon}) is far from reasonable when $p$ is large: the $\chi^2$ curve and the histogram, which are supposed to be matched, separate from each other with the increase of the value of $p.$ See the caption in Figure \ref{fig1} for more details.

The sizes and powers of the tests by using (\ref{afternoon}) and by Theorem \ref{yet} are estimated from our simulation and summarized in Table \ref{table_sphericity} at Section \ref{simulation_study}. A further analysis on this results is presented in the same section.

Finally, when $p \geq n$, we know the LRT does not exist as mentioned above. There are some recent works on choosing other statistics to study the spherical test of (\ref{have}), see, for example, Ledoit and Wolf (2002) and Chen, Zhang and Zhong (2010).

%\pagebreak
\subsection{Testing Independence of Components of Normal Distributions}\lbl{Independence}

Let $k\geq 2,\, p_1, \cdots, p_k$ be positive integers. Denote $p=p_1+\cdots + p_k$ and
\begin{eqnarray}
\bd{\Sigma}=(\bd{\Sigma}_{ij})_{p\times p}
\end{eqnarray}
be a positive definite matrix, where $\bd{\Sigma}_{ij}$ is a $p_i\times p_j$ sub-matrix for all $1\leq i,j \leq k.$ Let $N_p(\bd{\mu}, \bd{\Sigma})$ be a $p$-dimensional normal distribution. We are testing
\begin{eqnarray}\lbl{Chicago}
H_0: \bd{\Sigma}_{ij}=\bd{0}\ \mbox{for all}\ 1\leq i< j\leq k\ \ \ \mbox{vs}\ \ \ H_a: H_0\ \mbox{is not true}.
\end{eqnarray}
In other words, $H_0$ is equivalent to that $\bd{\xi}_1, \cdots, \bd{\xi}_k$ are independent, where $(\bd{\xi}_1', \cdots, \bd{\xi}_k')'$ has the distribution $N_p(\mu, \bd{\Sigma})$ and $\bd{\xi}_i\in \mathbb{R}^{p_i}$ for $1\leq i \leq k.$ Let $\bd{x}_1, \cdots, \bd{x}_N$ be i.i.d. with distribution $N_p(\mu, \bd{\Sigma}).$ Set $n=N-1.$ Let $\bd{S}$ be the covariance matrix as in (\ref{golden}). Now we partition $\bd{A}:=n\bd{S}$ in the following way:
\beaa
%\bar{\bd{X}}=\left(
%               \begin{array}{c}
%                 \bar{\bd{X}}_1\\
%                 \bar{\bd{X}}_2\\
%                 \vdots \\
%                 \bar{\bd{X}}_k\\
%               \end{array}
%             \right)
%\ \ \ \mbox{and}\ \ \
\bd{A}=
\left(
  \begin{array}{cccc}
    \bd{A}_{11} & \bd{A}_{12} & \cdots & \bd{A}_{1k} \\
    \bd{A}_{21} & \bd{A}_{22} & \cdots & \bd{A}_{2k} \\
    \vdots &  \cdots & \cdots  & \vdots \\
    \bd{A}_{k1} & \bd{A}_{k2} & \cdots & \bd{A}_{kk} \\
  \end{array}
\right)
\eeaa
where $\bd{A}_{ij}$ is a $p_i\times p_j$ matrix. Wilks (1935) shows that the likelihood ratio statistic for testing
(\ref{Chicago})  is given by
\bea\lbl{wired5}
\Lambda_n =\frac{|\bd{A}|^{(n+1)/2}}{\prod_{i=1}^k|\bd{A}_{ii}|^{(n+1)/2}}:= (W_n)^{(n+1)/2},
\eea
see also Theorem 11.2.1 from Muirhead (1982). Notice that $W_n=0$ if
 $p > N=n+1$, since the matrix $\bd{A}$ is not of full rank.  From (\ref{wired5}), we know that the LRT of level $\alpha$ for testing $H_0$ in (\ref{Chicago}) is $\{\Lambda_n\leq c_{\alpha}\}=\{W_n\leq c_{\alpha}'\}.$  Set
\begin{eqnarray}
f  =   \frac{1}{2} \Big( p^2 - \sum_{i=1}^k p_i^2 \Big)\ \ \mbox{and}\ \
\rho  =  1 - \frac{2 \Big( p^3 - \sum_{i=1}^k p_i^3 \Big) + 9 \Big( p^2 - \sum_{i=1}^k p_i^2 \Big)}{6(n+1) \Big( p^2 - \sum_{i=1}^k p_i^2 \Big)}. \nonumber
%\omega_2 & = & \frac{1}{\rho^2 n^2} \Bigg[ \frac{1}{48} \Big( p^4 - \sum_{i=1}^k p_i^4 \Big) - \frac{5}{96} \Big( p^2 - \sum_{i=1}^k p_i^2 \Big) - \frac{p^3-\sum_{i=1}^k p_i^3}{72 \big( p^2 - \sum_{i=1}^k p_i^2 \big) } \Bigg], \nonumber \\
%f & = &  \frac{1}{2} \Big( p^2 - \sum_{i=1}^k p_i^2 \Big), \ \ \ \mathrm{and} \ \ \ \gamma = (\rho n)^2 \omega_2. \nonumber
\end{eqnarray}
When $n$ goes to infinity while all $p_i$'s remain fixed, the traditional chi-square approximation to the  distribution of $\Lambda_n$ is referenced from Theorem 11.2.5 in Muirhead (1982):
\begin{eqnarray} \label{Independence_classic}
-2\rho \log \Lambda_n\ \mbox{converges to}\ \chi_f^2 \ \mbox{in distribution}
\end{eqnarray}
as $n\to\infty.$ Now we study the case when $p_i$'s are proportional to $n.$ For convenience of taking limit, we assume that $p_i$ depends on $n$ for each $1\leq i \leq k.$
\begin{theorem}\lbl{energy} Assume $n>p+1$ for all $n\geq 3$ and $p_i/n\to y_i\in (0, 1)$ as $n\to\infty$ for each $1\leq i\leq k.$ Recall $W_n$ as  defined in (\ref{wired5}). Then, under $H_0$ in (\ref{Chicago}),  $(\log W_n -\mu_n)/\sigma_n$ converges in distribution to $N(0,1)$ as $n\to\infty,$ where
\begin{eqnarray*}
 \mu_n=-r_{n-1}^2\Big(p-n+\frac{3}{2}\Big)+ \sum_{i=1}^k r_{n-1,i}^2\Big(p_i-n+\frac{3}{2}\Big)\  \mbox{and}\ \  \sigma_n^2=2r_{n-1}^2-2\sum_{i=1}^k r_{n-1,i}^2
\end{eqnarray*}
and  $r_x=(-\log (1-\frac{p}{x}))^{1/2}$ for $x>p$ and $r_{x,i}=(-\log (1-\frac{p_i}{x}))^{1/2}$ for $x>p_i$ and $1\leq i \leq k.$
\end{theorem}
%From the notation $p=p_1+\cdots + p_k$ given at the beginning of this section, we know that the assumptions in the above theorem %imply that $y:=\lim_{n\to\infty}p/(n-1)\in (0,1].$

%\underline{\textbf{Do we have to assume $n>1+p$ instead of $n>1+\max_{1\leq i \leq k}p_i$}}?\\
Though $H_0$ in (\ref{Chicago}) involves with unknown $\bd{\Sigma}_{ii}$'s, the limiting distribution in Theorem \ref{energy} is pivotal. This actually can be quickly seen by transforming $\bd{y}_i=\bd{\Sigma}^{-1/2}(\bd{x}_i-\mu)$ for $1\leq i \leq N.$ Then  $\bd{y}_1, \cdots, \bd{y}_N$ are i.i.d. with distribution $N_p(\bd{0}, \bd{I}_p).$ Put this into (\ref{wired5}), the $\bd{\Sigma}_{ii}$'s are then canceled in the fraction under the null hypothesis. See also the interpretation in terms of group transformations on p. 532 from Muirhead (1982).

\begin{figure}[t]
\centerline{\includegraphics[scale=0.44, bb=0 0 1000 525]{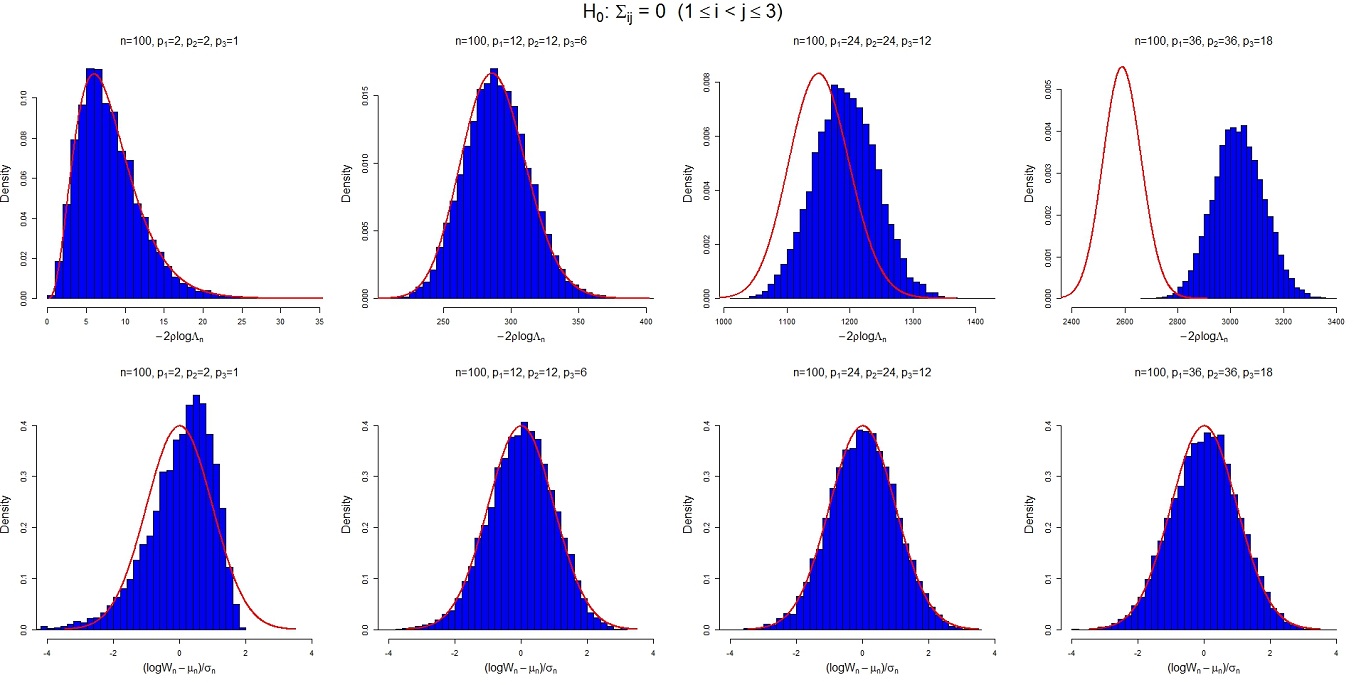}}
\caption{\sl Comparison between Theorem \ref{energy} and (\ref{Independence_classic}). We choose $k=3$, $n=100$ and $p=5, 30, 60, 90$ with $p_1:p_2:p_3=2:2:1$.  The pictures in the top row show that the histogram of $-2\rho\log \Lambda_n$ move away gradually from $\chi^2$ curve  as  $p$ grows. The pictures in the bottom row indicate that $(\log W_n-\mu_n)/\sigma_n$ and $N(0,1)$-curve match better as $p$ becomes larger. }
\label{fig2}
\end{figure}

We simulate the two cases in Figure \ref{fig2}: (i) the classical chi-square approximation (\ref{Independence_classic}); (ii) the central limit theorem based on Theorem \ref{energy}. The results show that when $p$ becomes large, the classical approximation in (\ref{Independence_classic}) is poor, however,  $(\log W_n -\mu_n)/\sigma_n$ in Theorem \ref{energy} fits the standard normal curve very well.

In Table \ref{table_independence_Component} from Section \ref{simulation_study}, we compare the sizes and powers of the two tests under the chosen $H_a$ explained in the caption. See the detailed explanations in the same section.

\subsection{Testing that Multiple Normal Distributions Are Identical}\lbl{Identical}

Given normal distributions $N_p(\mu_i, \bd{\Sigma}_i),\, i=1,2,\cdots, k,$ we are testing that they are all identical, that is,
\begin{eqnarray}\lbl{multiple}
H_0: \bd{\mu}_1=\cdots =\bd{\mu}_k,\ \bd{\Sigma}_1=\cdots = \bd{\Sigma}_k\ \ \mbox{vs}\ H_a: H_0\ \ \mbox{is not true}.
\end{eqnarray}
Let $\{\bd{y}_{ij};\, 1\leq i \leq k,\, 1\leq j \leq n_i\}$ be independent $p$-dimensional random vectors, and $\{\bd{y}_{ij};\, 1\leq j \leq n_i\}$ be i.i.d. from $N(\mu_i, \bd{\Sigma}_i)$ for each $i=1,2,\cdots, k.$ Set
\begin{eqnarray*}
& & \bd{A}=\sum_{i=1}^kn_i(\bar{\bd{y}}_i-\bar{\bd{y}})(\bar{\bd{y}}_i-\bar{\bd{y}})',\ \ \ \ \ \ \ \   \ \ \bd{B}_i=\sum_{j=1}^{n_i}(\bd{y}_{ij}-\bar{\bd{y}}_i)(\bd{y}_{ij}-\bar{\bd{y}}_i)'\ \ \mbox{and}\\
& & \bd{B}=\sum_{i=1}^k\bd{B}_i=\sum_{i=1}^k\sum_{j=1}^{n_i}(\bd{y}_{ij}-\bar{\bd{y}}_i)(\bd{y}_{ij}-\bar{\bd{y}}_i)'
\end{eqnarray*}
where
\begin{eqnarray*}
\bar{\bd{y}}_i=\frac{1}{n_i}\sum_{j=1}^{n_i}\bd{y}_{ij},\ \ \bar{\bd{y}}=\frac{1}{n}\sum_{i=1}^kn_i\bar{\bd{y}}_i,\ \ n=\sum_{i=1}^kn_i.
\end{eqnarray*}
%The following result is from Theorem 10.8.1 and Corollary 10.8.3 in \cite{Muirhead1982}.
The following likelihood ratio test statistic for  (\ref{multiple}) is first derived by Wilks (1932):
%According to Theorem 10.8.1 from \cite{Muirhead1982},  the LRT  statistic  for testing $H_0$ in (\ref{multiple}) is
\begin{eqnarray}\lbl{school}
\Lambda_n=\frac{\prod_{i=1}^k|\bd{B}_i|^{n_i/2}}{|\bd{A} + \bd{B}|^{n/2}}\cdot \frac{n^{pn/2}}{\prod_{i=1}^kn_i^{pn_i/2}}.
\end{eqnarray}
See also Theorem 10.8.1 from Muirhead (1982). The likelihood ratio test will reject the null hypothesis  if $\Lambda_n \leq c_{\alpha}$, where
the critical value $c_{\alpha}$ is determined so that the significance level of the test is equal
to $\alpha$. Note that when $p > n_i$, the matrix $\bd{B}_i$ is not of full rank for $i=1,\cdots, k,$
and consequently their determinants are equal to zero, so is the likelihood ratio
statistic $\Lambda_n$. Therefore,  to consider the test of (\ref{multiple}), one needs $p \leq  \min\{n_i;\, 1\leq i \leq k\}.$ Perlman (1980) shows that the LRT  is unbiased for testing $H_0.$  Let
\begin{eqnarray}
f = \frac{1}{2}p(k-1)(p+3) \ \ \ \mbox{and}\ \ \
\rho = 1- \frac{2p^2+9p+11}{6(k-1)(p+3)n}\Big(\sum_{i=1}^k\frac{n}{n_i} - 1 \Big).
\end{eqnarray}
When the dimension $p$ is considered fixed, the following asymptotic distribution of $\log \Lambda_n$ under the null hypothesis (\ref{multiple}) is a corollary from Theorem 10.8.4 in Muirhead (1982):
\begin{eqnarray} \label{classic_multiple_normal}
-2\rho \log \Lambda_n \ \ \mbox{converges to}\ \   \chi_f^2
\end{eqnarray}
in distribution as $\min_{1\leq i \leq k}{n_i}\to \infty.$ When $p$ grows with the same rate of $n_i$, we have the following theorem.
%In Section 10.3 from \cite{Anderson}, Anderson suggested the use of a modified test $\Lambda^*$ in %which the sample sizes $N_i$ are replaced with $n_i=N_i-1$ and $N$ is replaced with %$n=\sum_{i=1}^kn_i=N-p,$ that is,
%\bea\lbl{Lebron}
%\Lambda_N^*=\frac{\prod_{i=1}^k|\bd{B}_i|^{n_i/2}}{|\bd{A} + \bd{B}|^{n/2}}\cdot %\frac{n^{pN/2}}{\prod_{i=1}^kn_i^{pN_i/2}}.
%\eea
%For clarity of taking limit, we let $N_i$ be as a function of $p$ for all $1\leq i \leq k.$
\begin{theorem}\lbl{town} Let $n_i=n_i(p)>p+1$ for all $p\geq 1$ and $\lim_{p\to\infty}{p/n_i}=y_i\in (0, 1]$ for all $1\leq i \leq k.$ Let $\Lambda_n$ be as in (\ref{school}). Then, under $H_0$ in (\ref{multiple}),
\begin{eqnarray*}
\frac{\log \Lambda_n - \mu_n}{n\sigma_n}\ \ \mbox{converges in distribution to} \  N(0, 1)
\end{eqnarray*}
as $p\to\infty,$ where
\begin{eqnarray*}
& & \mu_n=\frac{1}{4}\Big[-2kp-\sum_{i=1}^ky_i +  nr_n^2(2p-2n+3)- \sum_{i=1}^k n_ir_{n_i'}^2(2p-2n_i+3)\Big],\\
%& = & \frac{1}{4}\Big[Nr_N^2(2p-2N+3)-\sum_{i=1}^k N_ir_{N_i'}^2(2p-2N_i+3)\Big],\ \ \ \ \ \ \ \ \ \ \ \ \ \ \ \\
& & \sigma_n^2=\frac{1}{2}\Big(\sum_{i=1}^k\frac{n_i^2}{n^2}r_{n_i'}^2-r_n^2\Big)>0,
\end{eqnarray*}
$n_i'=n_i-1$ and $r_x=\left(-\log \left(1-\frac{p}{x}\right)\right)^{1/2}$ for $x>p.$
\end{theorem}

The limiting distribution in Theorem \ref{town} is independent of $\mu_i$'s and $\bd{\Sigma}_i$'s. This can be seen by defining $\bd{z}_{ij}=\bd{\Sigma}_1^{-1/2}(\bd{y}_{ij}-\mu_1)$, we then know $\bd{z}_{ij}$'s are i.i.d. with distribution $N_p(\bd{0}, \bd{I}_p)$ under the null. It can be easily verified that the $\mu_i$'s and $\bd{\Sigma}_i$'s are canceled from the numerator and the denominator of    $\Lambda_n$ in (\ref{school}), and hence the right hand side only depends on $\bd{z}_{ij}$'s.

\begin{figure}[t]
\centerline{\includegraphics[scale=0.44, bb=0 0 1000 525]{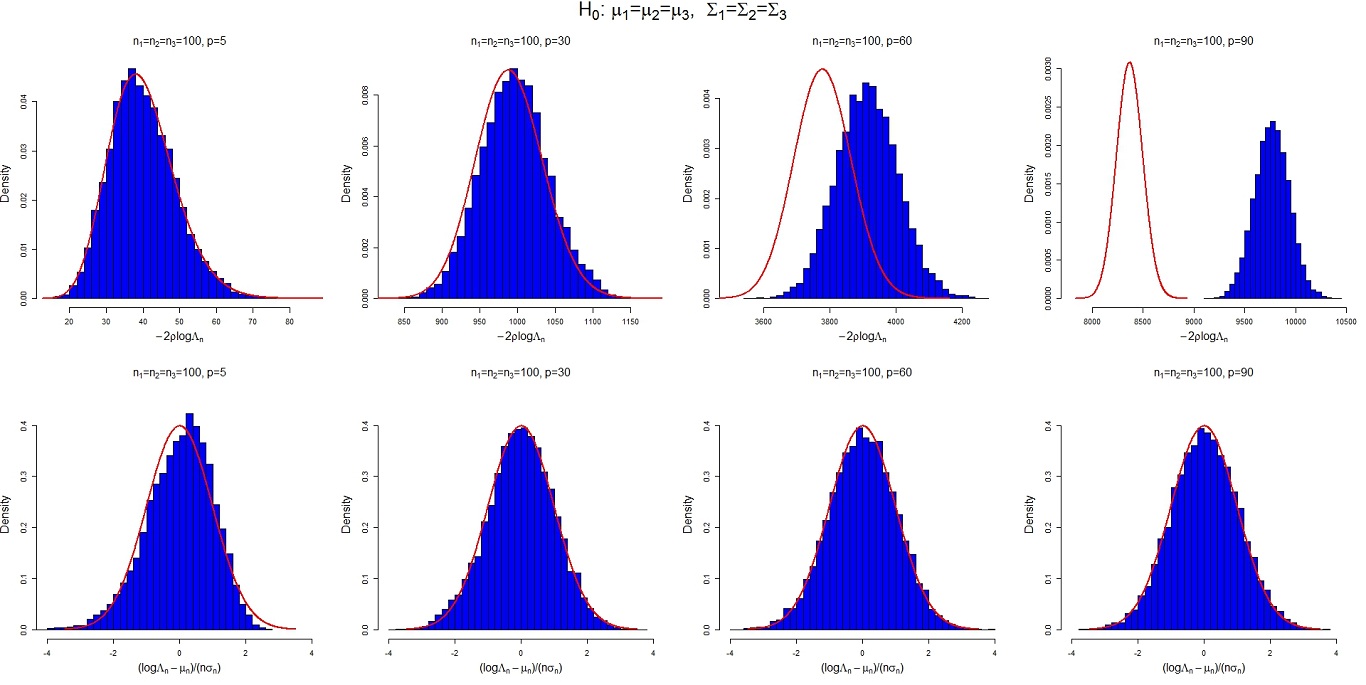}}
\caption{\sl Comparison between Theorem \ref{town} and (\ref{classic_multiple_normal}). We choose $n_1=n_2=n_3=100$ with $p=5, 30, 60, 90.$ The pictures in the top row show that the $\chi^2$ curves stay farther away  gradually from the histogram of $-2\rho\log \Lambda_n$  as $p$ grows. The pictures in the bottom row show that the $N(0, 1)$-curve fits the histogram of $(\log \Lambda_n-\mu_n)/(n\sigma_n)$  very well as $p$ becomes large.
%({\red change the value of $n$ in the pictures}).
}
\label{fig3}
\end{figure}

From the simulation shown in Figure \ref{fig3}, we see that when $p$ gets larger, the chi-square curve and the histogram are moving farther apart as $p$ becomes large, however, the normal approximation in Theorem \ref{town} becomes better.  The sizes and powers are estimated and summarized in Table \ref{table_Equality_3distribution} at Section \ref{simulation_study}. See more detailed explanations in the same section.

\subsection{Testing Equality of Several Covariance Matrices}\lbl{Constitution}
Let $k\geq 2$ be an integer. For $1\leq i \leq k,$ let $\bd{x}_{i1}, \cdots, \bd{x}_{in_i}$ be i.i.d. $N_p(\bd{\mu}_i, \bd{\Sigma}_i)$-distributed random vectors. We are considering
\begin{eqnarray}\lbl{bang}
H_0: \bd{\Sigma}_1=\cdots =\bd{\Sigma}_k\ \ \mbox{vs}\ \ \ H_a: H_0 \ \mbox{is not true}.
\end{eqnarray}
Denote
\begin{eqnarray*}
\overline{\bd{x}}_i=\frac{1}{n_i}\sum_{j=1}^{n_i}\bd{x}_{ij}\ \ \mbox{and}\ \ \ \bd{A}_i=\sum_{j=1}^{n_i}(\bd{x}_{ij}-\overline{\bd{x}}_i)(\bd{x}_{ij}-\overline{\bd{x}}_i)',\ \ \ \ 1\leq i \leq k,
\end{eqnarray*}
and
\begin{eqnarray*}
\bd{A}=\bd{A}_1 +\cdots + \bd{A}_k\ \ \ \mbox{and}\ \ \ \ n=n_1 +\cdots + n_k.
\end{eqnarray*}
Wilks (1932) gives the likelihood ratio test of (\ref{bang}) with a test statistic
\begin{eqnarray}\lbl{fellows}
\Lambda_n = \frac{\prod_{i=1}^k \left( \mathrm{det} \mathbf{A}_i \right)^{n_i/2}}{\left( \mathrm{det} \mathbf{A} \right)^{n/2}} \cdot \frac{n^{np/2}}{\prod_{i=1}^{k} n_i^{n_ip/2}}
\end{eqnarray} and the test rejects the null hypothesis $H_0$ at $\Lambda_n \leq c_\alpha$, where the critical value $c_\alpha$ is determined so that the test has the significance level of $\alpha$. Note that $\mathbf{A}_i$ does  not have a full rank when $p>n_i$ for any $i=1, \ldots, k,$ and hence their determinants are equal to zero. So the test statistic $\Lambda_n$ is not defined.  Therefore, we assume $p\leq n_i$ for all $i = 1, \ldots, k$ when study the likelihood ratio test of (\ref{bang}). The drawback of the likelihood ratio test is on its bias (see Section 8.2.2 of Muirhead (1982)). Bartlett (1937) suggests using a modified likelihood ratio test statistic $\Lambda_n^*$ by substituting every sample size $n_i$ with its degree of freedom $n_i-1$ and substituting the total sample size $n$ with $n-k$:
\begin{eqnarray} \label{stanford3}
\Lambda_n^* = \frac{\prod_{i=1}^k \left( \mathrm{det} \mathbf{A}_i \right)^{(n_i-1)/2}}{\left( \mathrm{det} \mathbf{A} \right)^{(n-k)/2}} \cdot \frac{(n-k)^{(n-k)p/2}}{\prod_{i=1}^{k} (n_i-1)^{(n_i-1)p/2}}.
\end{eqnarray}
The unbiased property of this modified likelihood ratio test is proved by Sugiura and Nagao (1968) for $k=2$ and by Perlman (1980) for a general $k$. Let
\begin{eqnarray}
f=\frac{1}{2}p(p+1)(k-1)\ \ \mbox{and}\ \ \rho = 1 - \frac{2p^2+3p-1}{6(p+1)(k-1)(n-k)} \Big( \sum_{i=1}^k \frac{n-k}{n_i-1} - 1 \Big). \nonumber
\end{eqnarray}
Box (1949) shows that when $p$ remains fixed, under the null hypothesis (\ref{bang}),
\bea\lbl{snake}
-2 \rho \log \Lambda_n^*\ \ \mbox{converges to}\ \ \chi^2_f
\eea
in distribution as $\min_{1\leq i \leq k}n_i\to \infty$ (See also Theorem 8.2.7 from Muirhead (1982)). Now, suppose $p$ changes with the sample sizes $n_i$'s. We have the following CLT.

%For convenience of taking limit, we regard $N_i$ as a function of $p$ for all $1\leq i \leq r.$
\begin{theorem}\lbl{eve} Assume $n_i=n_i(p)$ for all $1\leq i \leq k$ such that $\min_{1\leq i \leq k}n_i>p+1$ and $\lim_{p\to\infty} p/n_i=y_i\in (0, 1].$ Let $\Lambda_n^*$ be as in (\ref{stanford3}). Then, under $H_0$ in (\ref{bang}), $(\log\Lambda_n^* - \mu_n)/((n-k)\sigma_n)$ converges in distribution to $N(0, 1)$ as $p\to\infty,$
%\begin{eqnarray*}
%\frac{\log\Lambda_n^* - \mu_n}{n\sigma_n}\ \ \mbox{converges in distribution to}\ \ N(0, 1)
%\end{eqnarray*}
where
\begin{eqnarray*}
& & \mu_n=\frac{1}{4}\Big[(n-k)(2n-2p-2k-1)\log (1-\frac{p}{n-k}) \\
& & \ \ \ \ \ \ \ \ \ \ \ \ \ \ \ \ \ \ \ \ \ \ \ \ \ \ \ \ \ \ \ \ \ \ -\sum_{i=1}^k(n_i-1)(2n_i-2p- 3)\log (1-\frac{p}{n_i-1})\Big],\\
& & \sigma_n^2=\frac{1}{2}\Big[\log (1-\frac{p}{n-k})-\sum_{i=1}^{k}\Big(\frac{n_i-1}{n-k}\Big)^2\log (1-\frac{p}{n_i-1})\Big]>0.
\end{eqnarray*}
\end{theorem}
The limiting distribution in Theorem \ref{eve} is independent of $\mu_i$'s and $\bd{\Sigma}_i$'s. This is obvious: let   $\bd{y}_{ij}=\bd{\Sigma}_i^{-1/2}(\bd{x}_{ij}-\mu_i)$, then  $\bd{y}_{i}$'s are i.i.d. with distribution $N_p(\bd{0}, \bd{I}_p)$ under the null. From the cancelation of $\bd{\Sigma}_i$'s in $\Lambda^*_n$ from (\ref{stanford3}) we see that the distribution of $\Lambda^*_n$ is free of $\mu_i$'s and $\bd{\Sigma}_i$'s under $H_0.$

Bai et al. (2009) and Jiang et al. (2012) study Theorem \ref{eve}  for the case  $k=2.$ Theorem \ref{eve} generalizes their results for any $k\geq 2.$ Further, the first four authors impose the condition  $\max\{y_1, y_2\}<1$ which excludes the critical case $\max\{y_1, y_2\}=1.$ There is no such a restriction in Theorem \ref{eve}.

\begin{figure}[t]
\centerline{\includegraphics[scale=0.44, bb=0 0 1000 525]{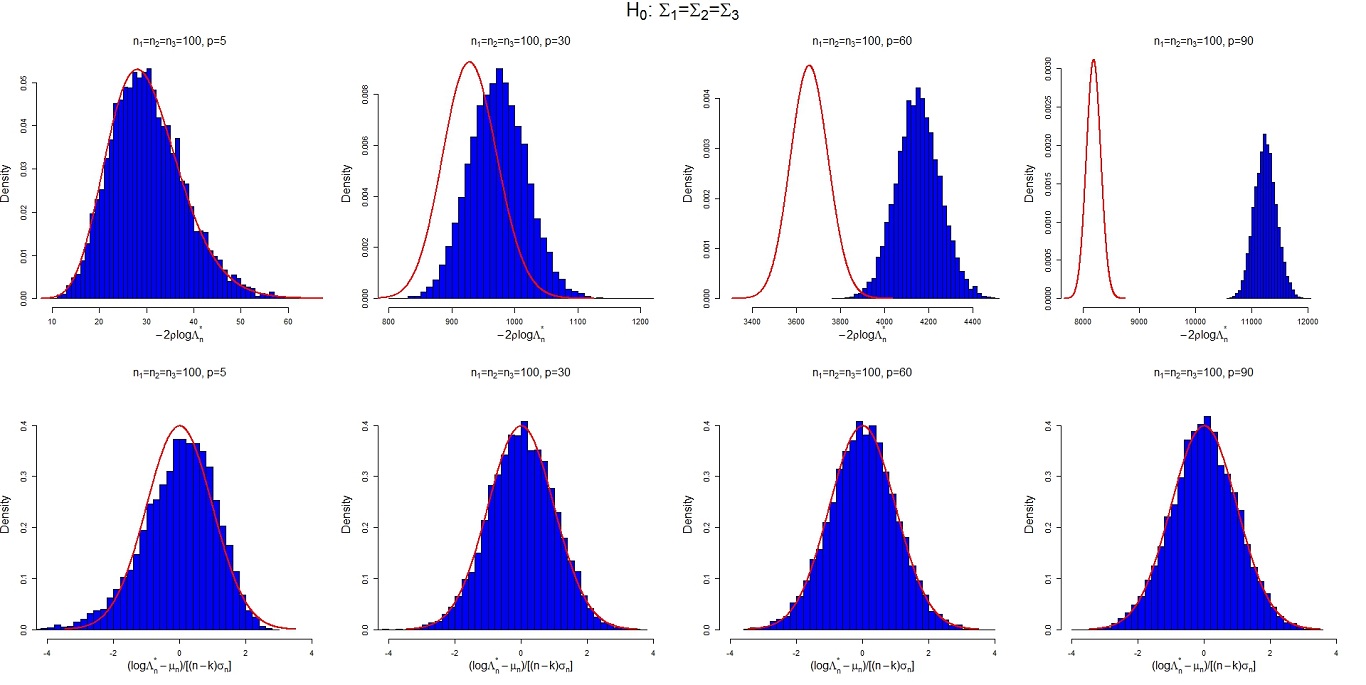}}
\caption{\sl Comparison between Theorem \ref{eve} and (\ref{snake}). We chose $n_1= n_2= n_3=100$ with $p=5, 30, 60, 90.$ The pictures in the top row show that the $\chi^2$ curves goes away quickly from the histogram of $-2\rho\log \Lambda_n^*$ as  $p$ grows. The pictures in the second row show that  the $N(0,1)$-curve fits the histogram of $(\log \Lambda^*_n-\mu_n)/[(n-k)\sigma_n]$  better as $p$ becomes larger.
%({\red change $n\sigma_n$ to $(n-k)\sigma_n$ and mention the values of $n_1, n_2, n_3$ in the pictures}).
}
\label{fig4}
\end{figure}

Figure \ref{fig4} presents our simulation with $k=3.$ It is interesting to see that the chi-square curve and the histogram almost separate from each other when $p$ is large, and at the same time the normal approximation in Theorem \ref{eve} becomes very good. In Table \ref{table_three_covariance} from Section \ref{simulation_study}, we estimate the sizes and powers of the two tests. The analysis is presented in the same section.

\subsection{Testing Specified Values for Mean Vector and Covariance Matrix}\lbl{Brown}

Let $\mathbf{x}_1, \cdots, \mathbf{x}_n$ be i.i.d. $\mathbb{R}^p$-valued random vectors from a normal distribution $N_p(\mathbf{\mu}, \mathbf{\Sigma})$, where $\mathbf{\mu} \in \mathbb{R}^p$ is the mean vector and $\mathbf{\Sigma}$ is the $p \times p$ covariance matrix. Consider the hypothesis test:
\begin{eqnarray*} \label{hypothesis_mean_covariance_specified_value_a}
H_{0}: & \mathbf{\mu}=\mathbf{\mu}_0 \ \mathrm{and} \ \mathbf{\Sigma}= \mathbf{\Sigma}_0
\ \ \mbox{vs}\ \ H_{a}: H_0\ \mbox{is not true},
\end{eqnarray*}
where $\mathbf{\mu}_0$ is a specified vector in $\mathbb{R}^p$ and $\mathbf{\Sigma}_0$ is a specified $p \times p$ non-singular matrix. By applying the transformation $\tilde{\mathbf{x}}_i = \mathbf{\Sigma}^{-1/2}(\mathbf{x}_i - \mathbf{\mu}_0)$, this hypothesis test is equivalent to the test of:
\begin{eqnarray} \label{hypothesis_mean_covariance_specified_value_b}
H_{0}: & \mathbf{\mu}=\mathbf{0} \ \mathrm{and} \ \mathbf{\Sigma}= \mathbf{I}_p
\ \ \mbox{vs}\ \ H_{a}: H_0\ \mbox{is not true}.
\end{eqnarray}
Recall the notation
\begin{eqnarray}\lbl{oar}
\mathbf{\bar{x}} = \frac{1}{n} \sum_{i=1}^{n} \mathbf{x}_i \ \ \mbox{and} \ \ \mathbf{A}=  \sum_{i=1}^{n}(\mathbf{x}_i-\mathbf{\bar{x}})(\mathbf{x}_i-\mathbf{\bar{x}})^{\prime}.
% \ \  \mbox{and}
%\ \ \mathbf{S}= \frac{\mathbf{A}}{n-1}.
\end{eqnarray}
The likelihood ratio test of size $\alpha$ of (\ref{hypothesis_mean_covariance_specified_value_b}) rejects $H_0$ if $\Lambda_n \leq c_{\alpha},$ where
%was studied by Anderson (see p.268 from Anderson (1958)) with the likelihood ratio statistic derived as
\begin{eqnarray}\label{likelihood_ratio_mean_covariance_specified_value}
\Lambda_n = \left( \frac{e}{n} \right)^{np/2} |\mathbf{A}|^{n/2} e^{ -\mathrm{tr}(\mathbf{A})/2} e^{-n\mathbf{\bar{x}}^\prime \mathbf{\bar{x}}/2}.
\end{eqnarray}
See, for example, Theorem 8.5.1 from Muirhead (1982). Note that the matrix $\mathbf{A}$ does not have a full rank when $p\geq n$ as discussed below (\ref{son}), therefore its determinant is equal to zero. This indicates that the likelihood ratio test of (\ref{hypothesis_mean_covariance_specified_value_b}) only exists when $p< n$. Sugiura and Nagao (1968) and Das Gupta (1969) show that this test with a rejection region $\{\Lambda_n \leq c_\alpha\}$ is unbiased, where the critical value $c_\alpha$ is chosen so that the test has the significance level of $\alpha$.

Theorem 8.5.5 from Muirhead (1982) implies that when the null hypothesis $H_0:\mathbf{\mu}=\mathbf{0}, \ \mathbf{\Sigma}= \mathbf{I}_p$ is true, the statistic
\bea\lbl{cargo}
-2 \rho \log \Lambda_n\ \ \mbox{converges to}\ \ \chi^2_f
\eea
as $n\to\infty$ with $p$ being fixed, where
\begin{eqnarray*}
\rho  =  1 - \frac{2p^2+9p+11}{6n(p+3)} \ \ \ \mbox{and} \ \ \ f = \frac{1}{2}p(p+3).
\end{eqnarray*}
Obviously, $\rho=\rho_n\to 1$ in this case. Davis (1971) improves the above result with a second order approximation. Nagarsenker and Pillai (1973) study the exact null distribution of $-2 \log \Lambda_n$ by using its moments.  Now we state our CLT result when $p$ grows with $n.$

%\noindent For clarity of taking limit, let $p=p_n$, that is, $p$ depends on $n.$
%\begin{theorem}\lbl{Ruth} Assume $n>p+1$ for all $n\geq 3$ and $\lim_{n\to\infty}\frac{p}{n}=y\in (0, 1].$ Then, under $H_0$ in (\ref{dong}),
%\begin{eqnarray*}
%\frac{\log Q_n - \mu_n}{n\sigma_n}\ \ \mbox{converges in distribution to}\ \ N(0, 1)
%\end{eqnarray*}
%as $n\to\infty,$ where
%\begin{eqnarray*}
%&  & \mu_n=-\frac{1}{4}\Big[n(2n-2p-3)\log (1-\frac{p}{n-1}) +2(n+1)p\Big]\ \ \mbox{and}\ \ \ \ \ \ \ \ \ \ \\
%&  &  \sigma_n^2= -\frac{1}{2}\Big(\frac{p}{n-1}+ \log(1-\frac{p}{n-1})\Big)>0.
%\end{eqnarray*}
%%with $r_n=(-\log (1-\frac{p}{n}))^{1/2}$ for all $n>p\geq 1.$
%\end{theorem}

\begin{figure}[t]
\centerline{\includegraphics[scale=0.44, bb=0 0 1000 525]{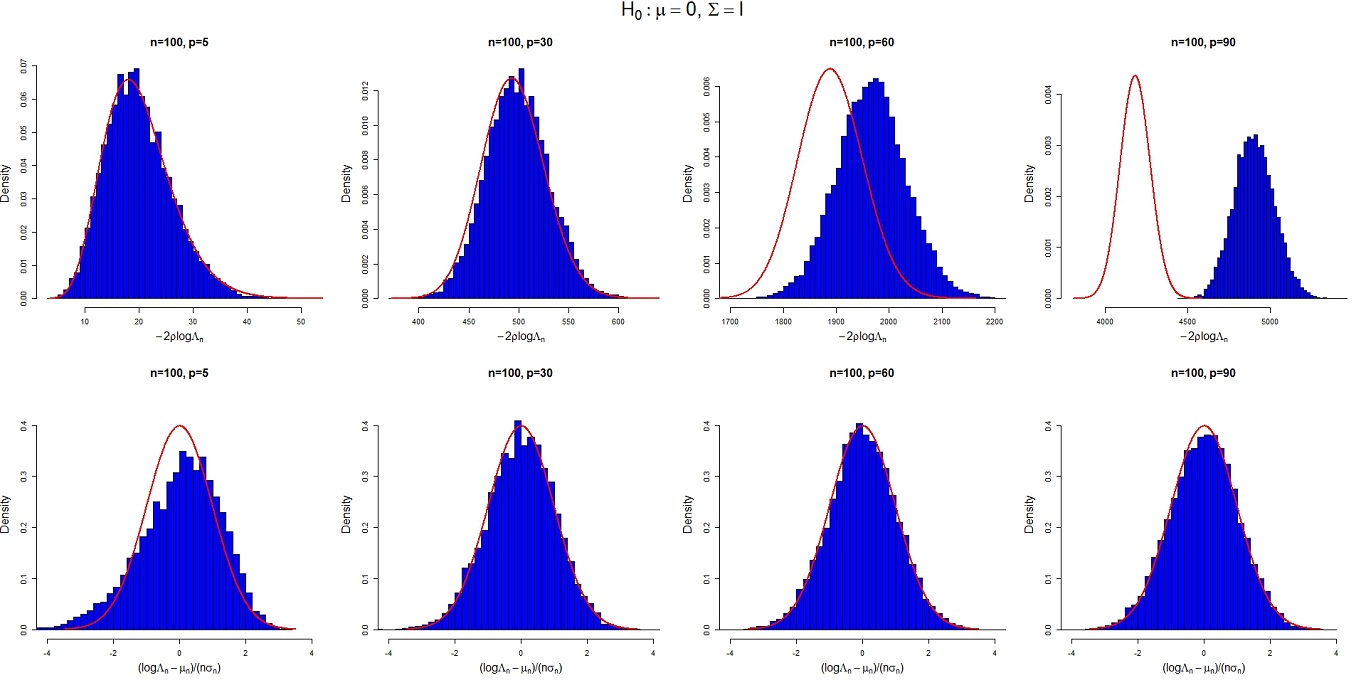}}
\caption{\sl Comparison between Theorem \ref{Ruth} and (\ref{cargo}). We chose $n=100$ with $p=5, 30, 60, 90.$ The pictures in the top row show that the $\chi^2$-curve stays away gradually from the histogram of $-2\log \Lambda_n$ as $p$ grows, whereas the $N(0, 1)$-curve fits  statistic $(\log \Lambda_n-\mu_n)/(n\sigma_n)$  better as shown from the bottom row. }
%{\red Change $Q_n$ to $\Lambda_n$.}
\label{fig5}
\end{figure}

\begin{theorem}\lbl{Ruth} Assume that $p:=p_n$  such that $n > 1+p$ for all $n \geq 3$ and $\lim_{n\to\infty}p/n = y\in (0, 1].$  Let $\Lambda_n$ be defined as in (\ref{likelihood_ratio_mean_covariance_specified_value}). Then under $H_{0}: \mathbf{\mu}=\mathbf{0}$ and $\mathbf{\Sigma}= \mathbf{I}_p$, $(\log \Lambda_n - \mu_n)/(n \sigma_n)$ converges in distribution to $N(0, 1)$ as $n\to\infty$,
where
\begin{eqnarray*}
&  & \mu_n=-\frac{1}{4}\Big[n(2n-2p-3)\log (1-\frac{p}{n-1}) +2(n+1)p\Big]\ \ \mbox{and}\ \ \ \ \ \ \ \ \ \ \\
&  &  \sigma_n^2= -\frac{1}{2}\Big(\frac{p}{n-1}+ \log(1-\frac{p}{n-1})\Big)>0.
\end{eqnarray*}
\end{theorem}
The simulations shown in Figure \ref{fig5} confirm that it is good to use Theorem \ref{Ruth} when $p$ is large and proportional to $n$ rather than the traditional chi-square approximation in (\ref{cargo}). In Table \ref{table_spicifed_distribution} from Section \ref{simulation_study}, we study the sizes and powers for the two tests based on  the $\chi^2$ approximation and our CLT. The understanding of the table is elaborated in the same section.

\subsection{Testing Complete Independence }\lbl{Correlation}

\noindent In this section, we study the likelihood ratio test  of the complete independence of the coordinates of a high-dimensional normal random vector. Precisely, let $\bd{R}=(r_{ij})_{p\times p}$ be the correlation matrix generated from $N_p(\bd{\mu}, \bd{\Sigma})$ and $\bd{x}=(x_1, \cdots, x_p)' \sim N_p(\bd{\mu}, \bd{\Sigma}).$ The test is
\bea\lbl{approach}
H_0: \bd{R}=\bd{I}\ \ \mbox{vs}\ \ H_a: \bd{R} \ne \bd{I}.
\eea
The null hypothesis $H_0$ is equivalent to that  $x_1, \cdots, x_p$ are independent or $\bd{\Sigma}$ is diagonal. To study the LRT, we need to understand the determinant of a sample correlation matrix generated by normal random vectors. In fact we will have a conclusion on the class of  spherical distributions, which is   more general than that of the normal distributions.   Let us first review two terminologies.

Let $\bd{x}=(x_1, \cdots, x_n)'\in \mathbb{R}^n$ and  $\bd{y}=(y_1, \cdots, y_n)'\in \mathbb{R}^n.$ Recall the {\it Pearson correlation coefficient} $r$ defined by
\bea\lbl{tube}
r=r_{\bd{x},\bd{y}}=\frac{\sum_{i=1}^n(x_i-\bar{\bd{x}})(y_i - \bar{\bd{y}})}
{\sqrt{\sum_{i=1}^n(x_i-\bar{\bd{x}})^2\cdot \sum_{i=1}^n(y_i-\bar{\bd{y}})^2}}
\eea
where $\bar{\bd{x}}=\frac{1}{n}\sum_{i=1}^n x_i$ and $\bar{\bd{y}}=\frac{1}{n}\sum_{i=1}^n y_i$.

We say a random vector $\bd{x}\in \mathbb{R}^n$ has a {\it spherical distribution} if $\bd{O}\bd{x}$ and $\bd{x}$ have the same probability distribution for all $n\times n$ orthogonal matrix $\bd{O}.$ Examples include the multivariate normal distribution $N_n(\bd{0}, \sigma^2\bd{I}_n)$, the ``$\epsilon$-contaminated" normal distribution $(1-\epsilon)N_n(\bd{0}, \bd{I}_n) + \epsilon N_n(\bd{0}, \sigma^2\bd{I}_n)$ with $\sigma>0$ and $\epsilon \in [0, 1]$, and the multivariate $t$ distributions.  See page 33 from Muirhead (1982) for more discussions.

Let $\bd{X}=(x_{ij})_{n\times p}=(\bd{x}_1, \cdots, \bd{x}_n)'=(\bd{y}_1, \cdots, \bd{y}_p)$ be an   $n\times p$ matrix such that $\bd{y}_1, \cdots,\bd{y}_p$ are independent random vectors with $n$-variate spherical distributions and $P(\bd{y}_i=0)=0$ for all $1\leq i \leq p$ (these distributions may be different). Let $r_{ij}=r_{\bd{y}_i, \bd{y}_j}$, that is, the Pearson correlation coefficient between $\bd{y}_i$ and $\bd{y}_j$ for $1\leq i
\leq j \leq n.$ Then,
\bea\lbl{playful}
\bd{R}_n :=(r_{ij})_{p\times p}
\eea
is the  sample correlation matrix. It is known that $\bd{R}_n$ can be written as $\bd{R}_n=\bd{U}'\bd{U}$ where $\bd{U}$ is an $n\times p$ matrix (see, for example, Jiang (2004a)). Thus, $\bd{R}_n$ does not have a full rank and hence $|\bd{R}_n|=0$ if $p>n.$ According to Theorem 5.1.3  from Muirhead (1982), the density function of $\bd{R}_n$ is given by
%proportional to  $|\bd{R}_n|^{(n-p-2)/2},\ \ (|r_{ij}|<1,\ i<j)$
\begin{eqnarray}\lbl{adopt}
\mbox{Constant} \cdot |\bd{R}_n|^{(n-p-2)/2}d\,\bd{R}_n.
%,\ \ (|r_{ij}|<1,\ i<j)
\end{eqnarray}
%with respect to the Lebesgue measure.

In the aspect of Random Matrix Theory, the limiting behavior of the largest eigenvalues of $\bd{R}_n$ and the empirical distributions of the eigenvalues of $\bd{R}_n$ are investigated by Jiang (2004a). For considering the construction of compressed sensing matrices, the statistical testing problems, the covariance structures of normal distributions, high dimensional regression in statistics and a wide range of applications including signal processing, medical imaging and seismology, the largest off-diagonal entries of $\bd{R}_n$ are studied by Jiang (2004b),  Li and Rosalsky (2006), Zhou (2007), Liu, Lin and Shao (2008), Li, Liu and Rosalsky (2009), Li, Qi and Rolsalski (2010) and Cai and Jiang (2011, 2012).

Let's now focus on the LRT of (\ref{approach}). According to p. 40 from Morrison (2004), the likelihood ratio test will reject the null hypothesis of (\ref{approach}) if
\begin{eqnarray}
|\mathbf{R}_n|^{n/2} \leq c_{\alpha}
\end{eqnarray} where $c_{\alpha}$ is determined so that the test has significant level of $\alpha$. It is also known (see, for example, Bartlett (1954) or p. 40 from Morrison (2005)) that when the dimension $p$ remains fixed and the sample size $n \to \infty$,
\begin{eqnarray} \label{kyookyoo}
-\Big(n-1-\frac{2p+5}{6}\Big) \log |\mathbf{R}_n| \overset{d}{\longrightarrow} \chi^2_{p(p-1)/2}.
\end{eqnarray}
This asymptotic result has been used for testing the complete independence of all the coordinates of a normal random vector in the traditional multivariate analysis when $p$ is small relative to $n$.

Now we study the LRT statistic when $p$ and $n$ are large and at the same scale. First, we give a general CLT result on spherical distributions.

\begin{theorem}\lbl{season} Let $p=p_n$ satisfy  $n\geq p+5$ and $\lim_{n\to\infty}p/n=y\in (0, 1].$ Let $\bd{X}=(\bd{y}_1, \cdots, \bd{y}_p)$ be an   $n\times p$ matrix such that $\bd{y}_1, \cdots,\bd{y}_p$ are independent random vectors with $n$-variate spherical distribution and $P(\bd{y}_i=\bd{0})=0$ for all $1\leq i \leq p$ (these distributions may be different). Recall $\bd{R}_n$ in (\ref{playful}). Then $(\log |\bd{R}_n| - \mu_n)/\sigma_n$ converges in distribution to  $N(0, 1)$
%\begin{eqnarray*}
%\frac{\log |\bd{R}_n| - \mu_n}{\sigma_n}\ \ \mbox{converges in distribution to}\  N(0, 1)
%\end{eqnarray*}
as $n\to\infty$, where
%$\mu_n=y-p-(p-n+1.5)r_n^2$ with $r_n=(-\log (1-\frac{p}{n-1}))^{-1/2}$ and
\begin{eqnarray*}
& &  \mu_n=(p-n+\frac{3}{2})\log (1-\frac{p}{n-1}) -\frac{n-2}{n-1}p\, ;\ \ \ \ \ \ \ \ \ \ \ \ \ \ \ \ \ \ \\
& & \sigma_n^2= -2\Big[\frac{p}{n-1} + \log \Big(1-\frac{p}{n-1}\Big)\Big].
\end{eqnarray*}
\end{theorem}
In the definition of $\sigma_n^2$ above, we need the condition $n\geq p+2.$
 %As mentioned earlier, $\bd{R}_n$ does not have a full rank if $p>n.$
However, the assumption ``$n\geq p+5$" still looks a bit stronger than ``$n\geq p+2$". In fact, we use the stronger one as a technical condition in the proof of Lemma \ref{cookup} which involves the complex analysis.

\begin{figure}[t]
\centerline{\includegraphics[scale=0.44, bb=0 0 1000 525]{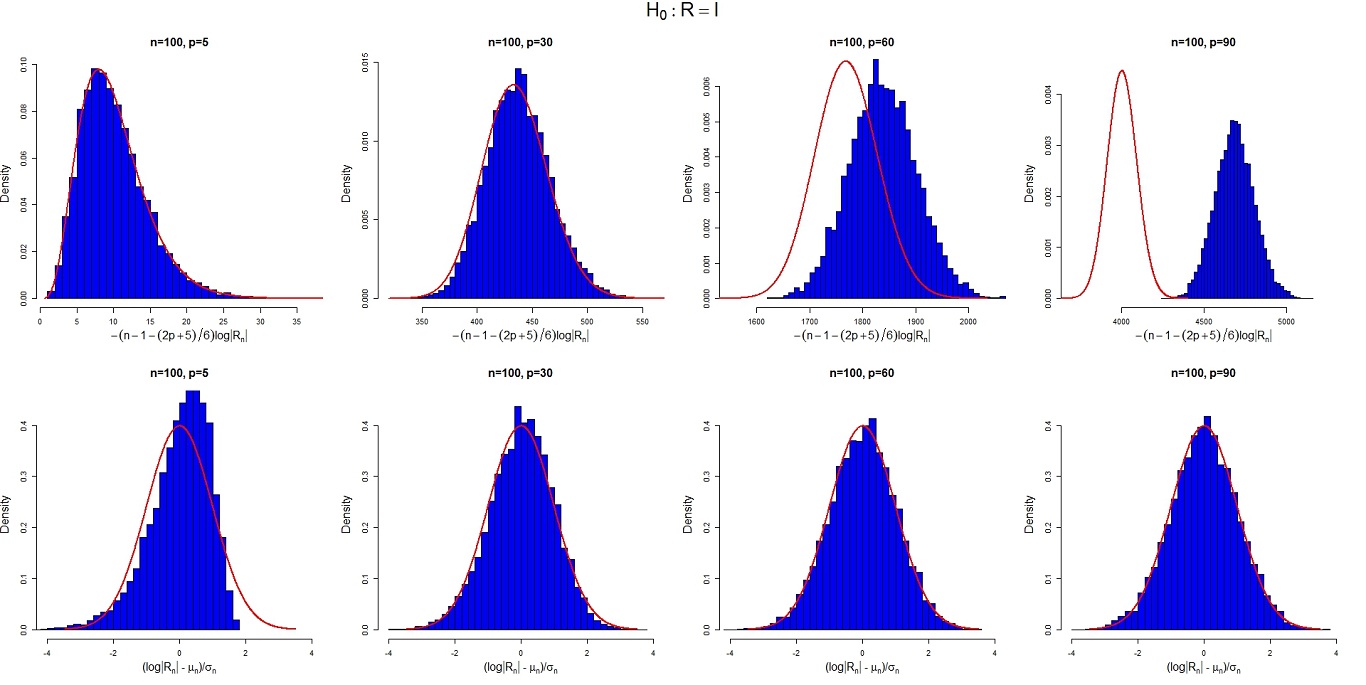}}
\caption{\sl Comparison between Corollary \ref{season_winter} and (\ref{kyookyoo}). We choose $n=100$ with $p=5, 30, 60, 90.$  The pictures in the first row show that, as $p$ becomes large, the $\chi^2$-curve fits the histogram of $-(n-1-\frac{2p+5}{6}) \log |\mathbf{R}_n|$ poorly. Those in the second row indicate that the N(0,1)-curve fits the histogram of $(\log |\mathbf{R}_n| - \mu_n)/\sigma_n$ very well as $p$ becomes large.}
\label{fig6}
\end{figure}

Notice that when the random vectors $\mathbf{x}_1, \ldots, \mathbf{x}_n$ are i.i.d. from a $p$-variate normal distribution $N_p(\mu, \mathbf{\Sigma})$ with complete independence (i.e., $\mathbf{\Sigma}$ is a diagonal matrix or the  correlation matrix  $\bd{R}=\mathbf{I}_p$). Write $\bd{X}=(x_{ij})_{n\times p}=(\bd{x}_1, \cdots, \bd{x}_n)'=(\bd{y}_1, \cdots, \bd{y}_p).$ Then, $\mathbf{y}_1, \cdots,\mathbf{y}_p$ are   independent random vectors from  $n$-variate normal distributions (these normal distributions may differ by their covariance matrices). It is also obvious that in this case $P(\mathbf{y}_i=0)=0$ for all $1\leq i \leq p$. Therefore, we have the following corollary.
\begin{coro} \label{season_winter}
Assume that $p:=p_n$ satisfy that $n\geq p+5$ and $\lim_{n\to\infty}p/n = y\in (0, 1].$  Let $\mathbf{x}_1, \cdots, \mathbf{x}_n$ be i.i.d. from $N_p(\mu, \mathbf{\Sigma})$ with the Pearson sample correlation matrix  $\mathbf{R}_n$ as defined in (\ref{playful}). Then, under $H_0$ in (\ref{approach}), $(\log |\mathbf{R}_n| - \mu_n)/\sigma_n$ converges in distribution to $N(0, 1)$ as $n\to\infty$, where
\begin{eqnarray*}
& &  \mu_n=\Big(p-n+\frac{3}{2}\Big)\log \left(1-\frac{p}{n-1}\right) -\frac{n-2}{n-1}p;\\
& & \sigma_n^2= -2\left[\frac{p}{n-1} + \log \left(1-\frac{p}{n-1}\right)\right] > 0.
\end{eqnarray*}
\end{coro}
According to Corollary \ref{season_winter}, the set $\{( \log |\mathbf{R}_n| - \mu_n)/\sigma_n \leq -z_{\alpha}\}$ is the rejection region with an asymptotic $1-\alpha$ confidence level for the LRT of (\ref{approach}), where the critical value $z_{\alpha}>0$ satisfies that $P(N(0,1)>z_{\alpha})=\alpha$ for all $\alpha \in (0,1).$ Figure \ref{fig6} shows that the chi-square approximation in (\ref{kyookyoo}) is good when $p$ is small, but behaves poorly as $p$ is large. At the same time, the normal approximation in Corollary \ref{season_winter} becomes better.

We simulate the sizes and powers of the two tests according to  the chi-square approximation in (\ref{kyookyoo}) and the CLT in Corollary \ref{season_winter} in Table \ref{table_complete_independence} at Section \ref{simulation_study}. See more analysis in the same section.

As mention earlier, when $p>n$, the LRT statistic $\log |\mathbf{R}_n|$ is not defined. So one has to choose other statistics rather than  $\log |\mathbf{R}_n|$ to study (\ref{approach}). See, for example, Schott (2005) and Cai and Ma (2012) for recent progress.
%In literature, Schott (2005) study the same test: let $T_{np} = (\sum_{i=2}^p \sum_{j=1}^{i-1} r_{ij}^2)- p(p-1)/(2n).$ He shows that if %$\lim_{n\to\infty}p/n = c \in (0, \infty)$, then  $T_{np}/\sigma_{np}$ converges in distribution to $N(0, 1)$ under $H_0$ in (\ref{approach}), %where $\sigma_{np} = p(p-1)(n-2)\cdot [(n-1)^2(n+1)]^{-1}.$  Since our main focus here is the size of the determinant of $\bd{R}_n,$ we leave the %comparison between our proposed test and Schott's one as a future work.

\section{Simulation Study: Sizes and Powers}\lbl{simulation_study}

In this part, for each of the six LRTs discussed earlier, we run simulation with 10,000 iterations to estimate the sizes and the powers of the LRTs using the CLT approximation and  the classical $\chi^2$ approximation. An analysis for each table is given. In the following discussion, the notation $\bd{J}_p$ stands for the $p\times p$ matrix whose entries are all equal to $1$ and $[x]$ stands for the integer part of $x>0.$

(1) {\it Table \ref{table_sphericity}}. This table corresponds to the sphericity  test that, for $N_p(\bd{\mu}, \bd{\Sigma})$, $H_0:\, \bd{\Sigma}=\lambda \bd{I}_p\ \ \mbox{vs}\ \ H_a:\, \bd{\Sigma} \ne \lambda \bd{I}_p$ with $\lambda$ unspecified. It is studied in Section \ref{Given_matrix}.  As expected, the $\chi^2$ approximation is good when $p$ is small relative to $n$, but not when $p$ is large. For example, at $n=100$ and $p=60,$ the size (type I error or alpha error)  for our normal approximation is $0.0511$ and power is $0.7914$, but the size for $\chi^2$ approximation is $0.3184$, which is too large to be used in practice. It is very interesting to see that our normal approximation is also as good as the $\chi^2$ approximation even when $p$ is small. Moreover, for $n=100$ and $p=90$ where the ratio $y=0.9$ is close to $1$,  the type I error in the CLT case is close to $5\%$ and the power is still decent at  $0.5406.$ Further, the power for the case of CLT drops as the ratio $p/n$ increases to $1$. This makes sense because the convergence rate of the CLT becomes slow. This can  be seen from Theorem \ref{yet} that $\sigma_n\to \infty$ as $p/n \to 1.$

(2) {\it Table \ref{table_independence_Component}}. In this table, we compare the sizes and powers of two tests under the chosen $H_a$ explained in the caption. The first  one is the classical $\chi^2$-approximation in (\ref{Independence_classic}) and the second is the CLT in Theorem \ref{energy} for the hypothesis that some  components of a normal distribution are independent. We observe  from the chart that our CLT approximation and the classical $\chi^2$ approximation are comparable for the small values of $p_i$'s. However, when $p_i$'s are large, noticing the last two rows in the table, our test is good whereas the $\chi^2$ approximation is no longer applicable because of the large sizes (type I errors). The power for the CLT drops when the values of $p_i$'s become large. This follows from  Theorem \ref{energy} that $\sigma_n\to\infty$ as $\sum p_i/n \to 1$, and hence the CLT-approximation does not perform well.

(3) {\it Table \ref{table_Equality_3distribution}}. We create this table for testing that several normal distributions are identical in Section \ref{Identical}. It is easily seen  that our CLT is good in all cases (except at the case of  $p=5$ where the type I error in our test is $0.0621$, slightly higher than $0.0512$ in the classical case). But when $p=60$ and $n_1=n_2=n_3=100$, the size in the classical case is $0.4542,$ too large to be used. It is worthwhile to notice that the power on the CLT becomes smaller as the value of $p$ becomes larger. This is easily understood from Theorem \ref{town} that the standard deviation diverges to infinity when $p/n \to 1.$ Equivalently, the convergence rate is poorer when $p$ gets closer to $n.$

(4) {\it Table \ref{table_three_covariance}}. This table relates to the test of the equality of the covariance matrices from $k$  normal distributions studied in  Section \ref{Constitution}. We take $k=3$ in our simulations. The sizes and powers of the chi-square approximation and the CLT in Theorem \ref{eve}  are summarized in the table. When $p=5$ and $n_1=n_2=n_3=100$, our CLT approximation gives a reasonable size of the test while the classical $\chi^2$ approximation is a bit better. However, for the same values of $n_i$'s, when $p=30,60,90,$  the size for the $\chi^2$ approximation is $0.2607, 0.9998$ and $1$, respectively, which are not recommended to be used in practice. Similar to the previous tests,  $\sigma_n \to \infty$ as $p/n\to 1$, where $\sigma_n^2$ is as in Theorem \ref{eve}. This implies that the convergence of the CLT is slow in this case. So it is not surprised to see that the power of the test based on the CLT in the table reduces as $p/n\to 1$.

(5) {\it Table \ref{table_spicifed_distribution}}. We generate this table by considering the LRT with $H_0: \bd{\mu}=\bd{0},\ \bd{\Sigma}=\bd{I}_p$ for the population distribution $N_p(\bd{\mu}, \bd{\Sigma}).$ The CLT is developed in Theorem \ref{Ruth}. In this table we study the sizes and powers for the two cases based on  the $\chi^2$ approximation and the CLT. At $n = 100, p = 5$ ($p$ is small), the $\chi^2$ test outperforms ours. The two cases are equally good at $n = 100, p = 30.$ When the values of $p$ are large at $60$ and $90,$ our CLT is still  good but the $\chi^2$ approximation is no longer useful. At the same time, it is easy to spot from the fourth column of the table that the power for the CLT-test drops as the ratio $p/n$ becomes large. It is obvious from Theorem \ref{Ruth} that the standard deviation $\sigma_n$ goes to infinity as the ratio approaches one. This causes the less precision when the sample size is not large.

(6) {\it Table \ref{table_complete_independence}}. This chart is created on the test that all of the components of a normal vector  are independent (but not necessarily identically distributed). It is studied in Corollary \ref{season_winter}. The sizes and powers of the two tests are estimated from simulation using the chi-square approximation in (\ref{kyookyoo}) and the CLT in Corollary \ref{season_winter}  from Section \ref{simulation_study} (the $H_a$ is explained in the caption). At all of the four cases of $n=100$ with $p=5,30,60$ and $90,$ the performance of our CLT-test is good, and it is even comparable with the classical $\chi^2$-test at the small value of $p=5.$ When $p=60$ and $90,$ the sizes of the $\chi^2$-test are  too big, while those of the CLT-test keep around $0.05.$ For the CLT-test itself, looking at the third and fourth rows of the table,  though the performance corresponding to $y=p/n=0.6$ is better than that corresponding to the high value of $y=p/n=0.9$ as expected, they are quite close. The only difference is the declining of the power as the rate $p/n$ increases. Again, this is easily seen from  Corollary \ref{season_winter} that the standard deviation $\sigma_n$ is divergent as $p$ is close to $n.$

\pagebreak

%\item In Section \ref{Conclusion_Discussion}, we summarize our results, state our method of the proofs and conclude by offering some open  problems.

%{\red The notation of $\bd{J}$ is a matrix whose entries are all equal to $1.$}

\begin{table}[!hbp]
\caption{Size and Power of LRT for Sphericity in Section \ref{Given_matrix}}

\begin{tabular}{|p{3cm}<{\centering}|p{2.5cm}<{\centering}|p{2.5cm}<{\centering}|p{2.5cm}<{\centering}|p{2.5cm}<{\centering}|}
  \hline
  & \multicolumn{2}{|c|}{Size under $H_0$} & \multicolumn{2}{|c|}{Power under $H_a$}  \\
  \cline{2-5}
  & CLT &  $\chi^2$ approx. & CLT &  $\chi^2$ approx.  \\
  \hline
  $n=100, p=5$ & 0.0562  & 0.0491 & 0.7525 & 0.7317   \\
  $n=100, p=30$ & 0.0581  & 0.0686 & 0.8700 & 0.8867   \\
  $n=100, p=60$ & 0.0511  & 0.3184 & 0.7914 & 0.9759   \\
  $n=100, p=90$ & 0.0518  & 1.0000 & 0.5406 & 1.0000   \\
  \hline
\end{tabular}
\label{table_sphericity}
\par
\noindent \footnotesize \sl The sizes (alpha errors) are estimated based on $10,000$ simulations from $N_p(\mathbf{0}, \mathbf{I}_p)$. The powers are estimated under the alternative hypothesis that $\mathbf{\Sigma} = diag(1.69,\cdots, 1.69, 1, \cdots, 1)$, where the number of 1.69 on the diagonal is equal to $[p/2]$.
\end{table}

\begin{table}[!hbp]
\caption{Size and Power of LRT for Independence of Three Components in Section \ref{Independence}}

\begin{tabular}{|p{5.8cm}<{\centering}|p{1.8cm}<{\centering}|p{1.8cm}<{\centering}|p{1.8cm}<{\centering}|p{1.8cm}<{\centering}|}
  \hline
   & \multicolumn{2}{|c|}{Size under $H_0$} & \multicolumn{2}{|c|}{Power under $H_a$}  \\
  \cline{2-5}
  & CLT &  $\chi^2$ approx. & CLT &  $\chi^2$ approx.  \\
  \hline
  $n=100, p_1=2, p_2=2, p_3=1$ & 0.0647  & 0.0458 & 0.7605 & 0.7176   \\
  $n=100, p_1=12, p_2=12, p_3=6$ & 0.0518  & 0.0543 & 0.9768 & 0.9778   \\
  $n=100, p_1=24, p_2=24, p_3=12$ & 0.0496  & 0.2171 & 0.8651 & 0.9757   \\
  $n=100, p_1=36, p_2=36, p_3=18$ & 0.0537  & 0.9998 & 0.4850 & 1.0000   \\
  \hline
\end{tabular}
\label{table_independence_Component}
\par
\noindent \footnotesize \sl The sizes (alpha errors) are estimated based on $10,000$ simulations from $N_p(\mathbf{0}, \mathbf{I}_p)$. The powers are estimated under the alternative hypothesis that $\mathbf{\Sigma} = 0.15\mathbf{J}_p + 0.85 \mathbf{I}_p$.
\end{table}

\begin{table}[!hbp]
\caption{Size and Power of LRT for Equality of Three  Distributions in Section \ref{Identical}}

\begin{tabular}{|p{5cm}<{\centering}|p{2.0cm}<{\centering}|p{2.0cm}<{\centering}|p{2.0cm}<{\centering}|p{2.0cm}<{\centering}|}
  \hline
   & \multicolumn{2}{|c|}{Size under $H_0$} & \multicolumn{2}{|c|}{Power under $H_a$}  \\
  \cline{2-5}
  & CLT &  $\chi^2$ approx. & CLT &  $\chi^2$ approx.  \\
  \hline
  $n_1=n_2=n_3=100, p=5$ & 0.0621  & 0.0512 & 0.7420 & 0.7135   \\
  $n_1=n_2=n_3=100, p=30$ & 0.0588  & 0.0743 & 0.8727 & 0.8936   \\
  $n_1=n_2=n_3=100, p=60$ & 0.0531  & 0.4542 & 0.6864 & 0.9770  \\
  $n_1=n_2=n_3=100, p=90$ & 0.0488  & 1.0000 & 0.3493 & 1.0000   \\
  \hline
\end{tabular}
\label{table_Equality_3distribution}
\par
\noindent \footnotesize \sl The sizes (alpha errors) are estimated based on $10,000$ simulations from three normal distributions of $N_p(\mathbf{0}, \mathbf{I}_p)$. The powers were estimated under the alternative hypothesis that $\mathbf{\mu}_1 = (0, \ldots, 0)^\prime$, $\mathbf{\Sigma}_1 = 0.5\mathbf{J}_p + 0.5 \mathbf{I}_p$; $\mathbf{\mu}_2 = (0.1, \ldots, 0.1)^\prime$, $\mathbf{\Sigma}_2 = 0.6\mathbf{J}_p + 0.4 \mathbf{I}_p$; $\mathbf{\mu}_3 = (0.1, \ldots, 0.1)^\prime$, $\mathbf{\Sigma}_3 = 0.5\mathbf{J}_p + 0.31 \mathbf{I}_p$.
\end{table}

\begin{table}[!hbp]
\caption{Size and Power of LRT for Equality of Three Covariance Matrices in Section \ref{Constitution}}

\begin{tabular}{|p{5cm}<{\centering}|p{2.0cm}<{\centering}|p{2.0cm}<{\centering}|p{2.0cm}<{\centering}|p{2.0cm}<{\centering}|}
  \hline
   & \multicolumn{2}{|c|}{Size under $H_0$} & \multicolumn{2}{|c|}{Power under $H_a$}  \\
  \cline{2-5}
  & CLT &  $\chi^2$ approx. & CLT &  $\chi^2$ approx.  \\
  \hline
  $n_1=n_2=n_3=100, p=5$ & 0.0805  & 0.0567 & 0.7157 & 0.6586   \\
  $n_1=n_2=n_3=100, p=30$ & 0.0516  & 0.2607 & 0.6789 & 0.9218   \\
  $n_1=n_2=n_3=100, p=60$ & 0.0525  & 0.9998 & 0.4493 & 1.0000  \\
  $n_1=n_2=n_3=100, p=90$ & 0.0535  & 1.0000 & 0.2297 & 1.0000   \\
  \hline
\end{tabular}
\label{table_three_covariance}
\par
\noindent \footnotesize \sl The sizes (alpha errors) are estimated based on $10,000$ simulations from $N_p(\mathbf{0}, \mathbf{I}_p)$. The powers are estimated under the alternative hypothesis that $\mathbf{\Sigma}_1 = \mathbf{I}_p$, $\mathbf{\Sigma}_2 = 1.1\mathbf{I}_p$, and $\mathbf{\Sigma}_3 = 0.9\mathbf{I}_p$.
\end{table}

\begin{table}[!hbp]
\caption{Size and Power of LRT for Specified Normal Distribution in Section \ref{Brown}}

\begin{tabular}{|p{3cm}<{\centering}|p{2.5cm}<{\centering}|p{2.5cm}<{\centering}|p{2.5cm}<{\centering}|p{2.5cm}<{\centering}|}
  \hline
  & \multicolumn{2}{|c|}{Size under $H_0$} & \multicolumn{2}{|c|}{Power under $H_a$}  \\
  \cline{2-5}
  & CLT &  $\chi^2$ approx. & CLT &  $\chi^2$ approx.  \\
  \hline
  $n=100, p=5$ & 0.0986  & 0.0471 & 0.5106 & 0.3818   \\
  $n=100, p=30$ & 0.0611  & 0.0657 & 0.7839 & 0.7898   \\
  $n=100, p=60$ & 0.0584  & 0.3423 & 0.7150 & 0.9583  \\
  $n=100, p=90$ & 0.0571  & 1.0000 & 0.4752 & 1.0000   \\
  \hline
\end{tabular}
\label{table_spicifed_distribution}
\par
\noindent \footnotesize \sl Sizes (alpha errors) are estimated based on $10,000$ simulations from $N_p(\mathbf{0}, \mathbf{I}_p)$. The powers are estimated under the alternative hypothesis that $\mathbf{\mu}=(0.1, \ldots, 0.1, 0, \ldots, 0)^\prime$ where the number of $0.1$ is equal to $[p/2]$ and $\mathbf{\Sigma} = \{\sigma_{ij}\}$ where $\sigma_{ij}=1$ for $i=j$, $\sigma_{ij}=0.1$ for $0<|i-j|\leq 3$, and $\sigma_{ij}=0$ for $|i-j|>3$.
\end{table}

\begin{table}[!hbp]
\caption{Size and Power of LRT for Complete Independence in Section \ref{Correlation}}

\begin{tabular}{|p{3cm}<{\centering}|p{2.5cm}<{\centering}|p{2.5cm}<{\centering}|p{2.5cm}<{\centering}|p{2.5cm}<{\centering}|}
  \hline
  & \multicolumn{2}{|c|}{Size under $H_0$} & \multicolumn{2}{|c|}{Power under $H_a$}  \\
  \cline{2-5}
  & CLT &  $\chi^2$ approx. & CLT &  $\chi^2$ approx.  \\
  \hline
  $n=100, p=5$ & 0.0548  & 0.0520 & 0.4311 & 0.4236   \\
  $n=100, p=30$ & 0.0526  & 0.0606 & 0.6658 & 0.6945   \\
  $n=100, p=60$ & 0.0522  & 0.3148 & 0.5828 & 0.9130  \\
  $n=100, p=90$ & 0.0560  & 1.0000 & 0.3811 & 1.0000   \\
  \hline
\end{tabular}
\label{table_complete_independence}
\par
\noindent \footnotesize \sl Sizes (alpha errors) are estimated based on $10,000$ simulations from $N_p(\mathbf{0}, \mathbf{I}_p)$. The powers are estimated under the alternative hypothesis that the correlation matrix $\mathbf{R}=(r_{ij})$ where $r_{ij}=1$ for $i=j$, $r_{ij}=0.1$ for $0<|i-j|\leq 3$, and $r_{ij}=0$ for $|i-j|>3$.
\end{table}

\pagebreak

\section{Conclusions and Discussions}\lbl{Conclusion_Discussion}
\noindent In this paper, we consider the likelihood ratio tests for the mean vectors and covariance matrices of high-dimensional normal distributions. Traditionally, these tests were performed by using the chi-square approximation. However, this approximation relies on a theoretical assumption that the sample size $n$ goes to infinity, while the dimension $p$ remains fixed.
%In practice, this requires the dataset to have a large sample size $n$ but a low dimension $p$.
As many modern datasets discussed in Section \ref{Introduction} feature high dimensions, these traditional likelihood ratio tests were shown to be less accurate in analyzing those datasets.

Motivated by the pioneer work of Bai et al. (2009) and Jiang et al. (2012), who prove two central limit theorems of the likelihood ratio test statistics for testing the high-dimensional covariance matrices of normal distributions, we examine in this paper other LRTs that are widely used in the multivariate analysis and prove the central limit theorems for their test statistics. By using the method developed in Jiang et al. (2012), that is, the asymptotic expansion of the multivariate Gamma function with high dimension $p$, we are able to derive the central limit theorems without relying on concrete random matrix models as demonstrated in Bai et al. (2009). Our method also has an advantage that the central limit theorems for the critical cases $\lim (p/n)=y=1$ or $\lim (p/n_i)=y_i=1$  are all derived, which is not the case in Bai et al. (2009) because of the restriction of their tools from the Random Matrix Theory. In real data analysis, as long as $n>p+1$ or $n_i>p+1$ in Theorems 1-5, or $n\geq p+5$ in Theorem \ref{season}, we simply take $y=p/n$ or $y_i=p/n_i$ to use the theorems. As Figures \ref{fig1}-\ref{fig6} and Tables 1-6 show, our CLT-approximations are all good even though  $p$ is relatively small.

The proofs in this paper are based on the analysis of the moments of the LRT statistics (five of six such moments are from literature and the last  one is derived by us as in  Lemma \ref{cookup}). The moment method we use here is different from that of the Random Matrix Theory employed in Bai et al. (2009) and  the Selberg integral used in Jiang et al. (2012).
%As mentioned earlier, when $y=\lim p/n$ or $y=\lim p/n_i$ is strictly larger  than $1$, the LRT test statistics do not exist simply because the %denominators in these statistics are functions of the determinants of certain matrices which are zero (the matrices are not of full ranks).

%However, since the proof by Bai et al. primarily uses the Random Matrix Theory that is shown less applicable to other LRTs with more complicated %test statistics, we develop a new method in our research which is based on the asymptotic expansion of the multivariate Gamma function with high %dimension $p$. It is shown that not only can our method simplify  the proof of Bai et al.'s result (see Jiang et al. (2012)), but it can also be %used to prove the central limit theorem for many other likelihood ratio test statistics.

Our research also brings out the following four interesting open problems:
\begin{enumerate}
\item All our central limit theorems in this paper are proved under the null hypothesis.
     As people want to assess the power of the test in many cases, it is also interesting to study the distribution of the test statistic under an alternative hypothesis. In the traditional case where $p$ is considered to be fixed while $n$ goes to infinity, the asymptotic distributions of many likelihood ratio statistics under the alternative hypotheses are derived by using the zonal polynomials (see,  e.g., Section 8.2.6, Section 8.3.4, Section 8.4.5 from Muirhead (1982)). It can be conjectured that in the high-dimensional case, there could be some new results regarding the limiting  distributions of the test statistics under the  alternative hypotheses. However, this is non-trivial and may require more investigation of the high-dimensional zonal polynomials. Some new understanding about the connection between the random matrix theory and the Jack polynomials (the zonal polynomials, the Schur polynomials and the zonal spherical functions are special cases) is  given by Jiang and Matsumoto (2011). A recent work by  Bai et al. (2009) study the high-dimensional LRTs through the random matrix theory. So the connection among the random matrix theory, LRTs and the Jack polynomials is obvious.  We are almost sure that the understanding by Jiang and Matsumoto (2011) will be useful in exploring the LRT  statistics under the  alternative hypotheses.

\item Except Theorem \ref{season} where the condition $n \geq  p+5$ is imposed due to a technical constraint, all other five central limit theorems in this paper are proved under the condition $n>p+1$ or $n_i> p+1$. This is because when this is not the case in the five theorems,  the likelihood ratio statistics will become undefined in these five cases. This indicates that tests other than the likelihood ratio ones shall be developed for analyzing a dataset with $p$ greater than $n$. For recent progress, see, for example,  Ledoit and Wolf (2002) and Chen et al. (2010) for the sphericity test, and Schott (2001, 2007) for testing the equality of multiple covariance matrices and Srivastava (2005) for testing the covariance matrix of a normal distribution. A power study for sphericity test is tried by  Onatski et al.  Despite these enlightening work mentioned above, other hypothesis tests for $p > n$ or $p>n_i$ are still an open area with many interesting problems to be solved.

\item In this paper we consider the cases when $p$ and $n$ or $n_i$ are proportional to each other, that is, $\lim p/n= y\in (0,1]$ or $\lim p/n_i= y_i\in (0,1].$ In practice, $p$ may be large but may not be large enough to be at the same scale of $n$ or $n_i$. So it is useful to derive the central limit theorems appeared in this paper under the assumption that $p\to\infty$ such that $p/n\to 0$ or $p/n_i\to 0$.

\item To understand the robustness of the six likelihood tests in this paper, one has to study the limiting behaviors of the LRT  statistics without the normality assumptions. This is feasible. For example, in Section \ref{Independence} we test the independence of several components of a normal distribution. The LRT statistic $W_n$ in (\ref{wired5}) can be written as the product of some independent random variables, say, $V_i$'s with beta distributions (see, e.g., Theorem 11.2.4 from Muirhead (1982)).  Therefore, it is possible that we can derive the CLT of $W_n$ for general $V_i$'s with the same means and variances as those of the beta distributions.

    %to test  $H_0: \mu_1=\cdots=\mu_k, \Sigma_1=\cdots=\Sigma_k$, it is well known that, under normality assumptions, the %likelihood ratio statistic $V_1$ or $W$ can be expressed as the product of several independent Beta variables, see Theorems %10.4.2 and 10.4.3 on Page 420-421 of Anderson (2003). Therefore, when the dimension $p$ and sample size $n$ increases %proportionally, the CLT's of $V_1$ or $W$ can easily follow by applying the delta method to the logarithms of the independent %products of Beta variables. (Note: $\beta(m,n)$ can be expressed as $\gamma(m)/(\gamma(m)+\gamma(n))$ and $\gamma(n)$ can be %expressed as a sum of the squares of $n$ independent standard normal variables) I believe that it may be of more interests to %establish the same results without normality assumptions.
\end{enumerate}

Finally, it is worthwhile to mention that some  recent works consider similar problems under the nonparametric setting, see, e.g., Cai et al. (2013), Cai and Ma, Chen et al. (2010), Li and Chen (2012), Qiu and Chen (2012) and Xiao and Wu (2013).

\section{Proofs}\lbl{Proofs}
This section is divided into some subsections. In each of them we prove a theorem introduced in Section \ref{Introduction}. We first develop some tools. The following are some standard notation.

For two sequence of numbers $\{a_n;\, n\geq 1\}$ and $\{b_n;\, n\geq 1\}$, the notation $a_n=O(b_n)$ as $n\to\infty$ means that $\limsup_{n\to\infty}|a_n/b_n|<\infty.$ The notation $a_n=o(b_n)$ as $n\to\infty$ means that $\lim_{n\to\infty}a_n/b_n=0.$  For two functions $f(x)$ and $g(x)$, the notation  $f(x)=O(g(x))$ and $f(x)=o(g(x))$ as $x\to x_0\in [-\infty, \infty]$ are similarly interpreted.

Throughout the paper $\Gamma(z)$ is the Gamma function defined on the complex plane $\mathbb{C}.$

\subsection{A Preparation}\lbl{preparation}
\setcounter{equation}{0}

\begin{lemma}\lbl{stirling} Let $b:=b(x)$ be a real-valued function defined on $(0, \infty)$.  Then,
%\begin{eqnarray*}
%\log \frac{\Gamma(x+b)}{\Gamma(x)} =b\log x+ \frac{b^2-b}{2x} + O\left(\frac{1}{\sqrt{x}}\right)
%\end{eqnarray*}
\begin{eqnarray*}
\log \frac{\Gamma(x+b)}{\Gamma(x)} =b\log x+ \frac{b^2-b}{2x} + c(x)
% O\left(\frac{1}{\sqrt{x}}\right)
\end{eqnarray*}
as $x \to +\infty$, where
\beaa
c(x)=\begin{cases}
O(x^{-1/2}), & \text{if\, $b(x)=O(\sqrt{x}\,)$;}\\
O(x^{-2}), & \text{if\, $b(x)=O(1).$}
\end{cases}
\eeaa
Further, for any constants $d>c,$ as $x\to+\infty$,
\beaa
\sup_{c\leq t \leq d}\Big|\log \frac{\Gamma(x+t)}{\Gamma(x)} -t\log x\Big| \to 0.
\eeaa
\end{lemma}
\textbf{Proof}. Recall the Stirling formula (see, e.g., p. 368 from Gamelin (2001) or (37) on
p. 204 from Ahlfors (1979)):
\begin{eqnarray}\lbl{lastflavor}
\log\Gamma(x)=\big(x-\frac{1}{2}\big)\, \mbox{log} x - x + \log \sqrt{2\pi}
+\frac{1}{12x} +O\left(\frac{1}{x^3}\right)
\end{eqnarray}
as $x\to +\infty$. We have that
%, where $\mbox{Log}\, z = \log |z| + i\theta$ with $\mbox{Arg}(z)=\theta\in (-\pi, \pi].$ Since $\mbox{Arg}(x)=0$ for $x>0$, we have
\begin{eqnarray}
& & \log \frac{\Gamma(x+b)}{\Gamma(x)}=  (x+b)\log (x+b) -x\log x-b-\frac{1}{2}\left(\log (x+b)-\log x\right) \nonumber\\
& & \ \ \ \ \ \ \ \ \ \ \ \ \ \ \ \ \ \ \ \ \   \ \ \ \ \ \ \ \ \ \ \ \ \ \ \ \ \ \ \ \ \ \ + \frac{1}{12}\left(\frac{1}{x+b}-\frac{1}{x}\right)+O\Big(\frac{1}{x^3}\Big)\ \ \ \ \ \ \lbl{alpha}
\end{eqnarray}
 as $x\to+\infty$. First, use the fact that $\log(1+t)= t -(t^2/2) +O(t^3)$ as $t\to 0$ to get
\begin{eqnarray*}
(x+b)\log (x+b) -x\log x &= & (x+b)\left(\log x +\log \Big(1+\frac{b}{x}\Big)\right) -x\log x\\
& = & (x+b)\Big(\log x +\frac{b}{x} -\frac{b^2}{2x^2} +O\Big(\frac{b^3}{x^3}\Big)\Big)-x\log x\\
& = & b\log x + b +\frac{b^2}{2x} +O\Big(\frac{b^3}{x^2}\Big) -\frac{b^3}{2x^2} + O\Big(\frac{b^4}{x^3}\Big)\\
& = & b\log x + b + \frac{b^2}{2x}+ c_1(x)
\end{eqnarray*}
as $x\to+\infty,$ where
\beaa
c_1(x)=
\begin{cases}
O(x^{-1/2}) & \text{if $b(x)=O(\sqrt{x})$;}\\
O(x^{-2}) & \text{if $b(x)=O(1)$.}
\end{cases}
\eeaa
%as $x\to +\infty.$
Similarly, as $x\to +\infty,$
\begin{eqnarray*}
& & \log (x+b)-\log x=\log \Big(1+\frac{b}{x}\Big)=
\begin{cases}
\frac{b}{x} + O(x^{-1}) & \text{if $b(x)=O(\sqrt{x})$;}\\
\frac{b}{x} + O(x^{-2}) & \text{if $b(x)=O(1)$}
\end{cases}
\end{eqnarray*}
and
\begin{eqnarray*}
 \frac{1}{x+b}-\frac{1}{x}=-\frac{b}{x(x+b)}=
\begin{cases}
O(x^{-3/2}) & \text{if $b(x)=O(\sqrt{x})$;}\\
O(x^{-2}) & \text{if $b(x)=O(1)$.}
\end{cases}
\end{eqnarray*}
%as $x\to+\infty.$
Substituting  these two assertions in (\ref{alpha}), we have
\begin{eqnarray}\lbl{publisher}
\log \frac{\Gamma(x+b)}{\Gamma(x)} =b\log x+ \frac{b^2-b}{2x} + c(x)
\end{eqnarray}
with
\beaa
c(x)=
\begin{cases}
O(x^{-1/2}) & \text{if $b(x)=O(\sqrt{x})$;}\\
O(x^{-2}) & \text{if $b(x)=O(1)$}
\end{cases}
\eeaa
as $x\to +\infty.$

For the last part, reviewing the whole proof above, we have  from (\ref{publisher}) that
\beaa
\log \frac{\Gamma(x+t)}{\Gamma(x)} =t\log x+ \frac{t^2-t}{2x} + O(x^{-2})
\eeaa
as $x\to +\infty$ uniformly for all $c\leq t \leq d.$ This implies the conclusion.\ \ \ \ \ \ \ \ $\blacksquare$

\begin{lemma}\lbl{software} Given $a>0.$ Define
\beaa
\eta(t)=\sup_{x\geq a}\Big|\log \frac{\Gamma(x+t)}{\Gamma(x)} -t\log x\Big|
\eeaa
for all $t> -a.$ Then $\lim_{t\to 0}\eta(t)=0.$
\end{lemma}
\textbf{Proof}. Let $d>c>0$ be two constants. Since $\Gamma(x)>0$ is continuous on $(0, \infty),$ then $g(x):=\log \Gamma(x)$ is uniformly continuous over the compact interval $[c/2, 2d].$ It then follows that
\bea\lbl{perron}
\sup_{c\leq x \leq d}\Big|\log \frac{\Gamma(x+\epsilon)}{\Gamma(x)}\Big| =\sup_{c\leq x \leq d}\Big|g(x+\epsilon)-g(x)\Big|\to 0
\eea
as $\epsilon \to 0.$ On the other hand, by the second part of Lemma \ref{stirling}, for any $\epsilon>0$, there exists $x_0>2a$ such that
\beaa
\sup_{|t|\leq a}\Big|\log \frac{\Gamma(x+t)}{\Gamma(x)} -t\log x\Big| < \epsilon
\eeaa
for all $x\geq x_0.$ Therefore,
\beaa
\sup_{x\geq x_0}\sup_{|t|\leq a}\Big|\log \frac{\Gamma(x+t)}{\Gamma(x)} -t\log x\Big| \leq \epsilon.
\eeaa
Then,
\beaa
\eta(t) &\leq & \epsilon +\sup_{a\leq x \leq x_0}\Big|\log \frac{\Gamma(x+t)}{\Gamma(x)} -t\log x\Big|\\
& \leq & \epsilon + (|\log a| + |\log x_0|)\cdot |t| + \sup_{a\leq x \leq x_0}\Big|\log \frac{\Gamma(x+t)}{\Gamma(x)}\Big|
\eeaa
for all $|t|\leq a.$ Consequently, we have from (\ref{perron}) that $\limsup_{t\to 0}\eta(t)\leq \epsilon$ for all $\epsilon>0,$ which concludes the lemma.\ \ \ \ \ \ \ $\blacksquare$

\begin{prop}(Proposition 2.1 from Jiang et al. (2012))\lbl{pizza} Let $n>p=p_n$ and $r_n=(-\log (1-\frac{p}{n}))^{1/2}.$ Assume that $p/n\to y\in (0, 1]$ and $t=t_n=O(1/r_n)$ as $n\to\infty.$  Then, as $n\to\infty,$
\begin{eqnarray*}
 \log \prod_{i=n-p}^{n-1}\frac{\Gamma(\frac{i}{2}-t)}{\Gamma(\frac{i}{2})}= pt(1+\log 2-\log n) + r_n^2\Big(t^2 +(p-n+1.5)t\Big) + o(1).
\end{eqnarray*}
\end{prop}

\begin{lemma}\lbl{card} Let  $n>p=p_n$ and $r_n=(-\log (1-\frac{p}{n}))^{1/2}.$ Assume $\frac{p}{n}\to y\in (0, 1]$ and $t=t_n=O(1/r_n)$ as $n\to\infty.$ Then
\begin{eqnarray}\lbl{cough}
\log \Big[ \frac{\Gamma(\frac{n}{2}+t)}{\Gamma(\frac{n}{2})}\cdot \frac{\Gamma(\frac{n-p}{2})}{\Gamma(\frac{n-p}{2} +t)}\Big] =r_n^2t + o(1)
\end{eqnarray}
as $n\to\infty.$
\end{lemma}
\textbf{Proof}. We prove the lemma by considering two cases.

\noindent{\it Case (i):  $y \in (0,1)$.} In this case, $n-p\to\infty$ and $\lim_{n\to\infty}r_n=(-\log (1-y))^{1/2}\in (0, \infty)$, and hence $\{t_n\}$ is bounded. By Lemma \ref{stirling},
\beaa
& & \log \frac{\Gamma(\frac{n}{2}+t)}{\Gamma(\frac{n}{2})}=t\log \frac{n}{2} +O\Big(\frac{1}{n}\Big) \\
& & \log \frac{\Gamma(\frac{n-p}{2})}{\Gamma(\frac{n-p}{2} +t)}=-t\log \frac{n-p}{2}+O\Big(\frac{1}{n-p}\Big)
\eeaa
as $n\to\infty.$ Add the two assertions up, we get that the left hand side of (\ref{cough}) is equal to
\bea\lbl{comfort}
-t\log \big(1-\frac{p}{n}\Big) + o(1)= r_n^2t + o(1)
\eea
as $n\to \infty.$ So the lemma holds for $y \in (0,1).$

\noindent{\it Case (ii):  $y =1$.} In this case, $r_n\to +\infty$ and $t_n\to 0$ as $n\to \infty.$ Recalling Lemma \ref{software}, we know that
\beaa
\Big|\log  \frac{\Gamma(\frac{n-p}{2}+t_n)}{\Gamma(\frac{n-p}{2} )} - t_n\log \frac{n-p}{2}\Big| \leq \eta(t_n)\to 0
\eeaa
as $n\to\infty$ by taking $a=1/2$ since $n-p\geq 1.$ That is,
\bea\lbl{data}
\log  \frac{\Gamma(\frac{n-p}{2} )}{\Gamma(\frac{n-p}{2}+t_n)} = -t_n\log \frac{n-p}{2} + o(1)
\eea
as $n\to\infty.$ By Lemma \ref{stirling} and the fact that $\lim_{n\to\infty}t_n=0$,
\beaa
\log \frac{\Gamma(\frac{n}{2}+t_n)}{\Gamma(\frac{n}{2})}=t_n\log \frac{n}{2} +o(1)
\eeaa
as $n\to\infty.$ Adding up the above two terms, then using the same argument as in (\ref{comfort}), we obtain (\ref{cough}).\ \ \ \ \ \ \ \ $\blacksquare$\\

\noindent Define
\begin{eqnarray}\lbl{complex}
\Gamma_p(z):=\pi^{p(p-1)/4}\prod_{i=1}^{p}\Gamma\Big(z-\frac{1}{2}(i-1)\Big)
\end{eqnarray}
for complex number $z$ with $\mbox{Re}(z)>\frac{1}{2}(p-1).$ See p. 62 from Muirhead (1982).

\begin{lemma}\lbl{by} Let $\Gamma_p(z)$ be as in (\ref{complex}). Let  $n>p=p_n$ and $r_n=(-\log (1-\frac{p}{n}))^{1/2}.$ Assume $\frac{p}{n}\to y\in (0, 1]$ and $s=s_n=O(1/r_n)$ and $t=t_n=O(1/r_n)$ as $n\to\infty.$ Then
\begin{eqnarray*}
\log \frac{\Gamma_p(\frac{n}{2} + t)}{\Gamma_p(\frac{n}{2}+s)}
 = p(t-s)(\log n -1-\log 2) + r_n^2\Big[(t^2-s^2) -\Big(p-n+\frac{1}{2}\Big)(t-s)\Big] +o(1)
\end{eqnarray*}
as $n\to\infty.$
\end{lemma}
\textbf{Proof}. First,
\beaa
\Gamma_p\Big(\frac{n}{2} + t\Big) & = & \pi^{p(p-1)/4}\prod_{i=1}^{p}\Gamma\Big(\frac{n-i}{2}+ t +\frac{1}{2}\Big)\\
& = & \pi^{p(p-1)/4}\prod_{j=n-p}^{n-1}\Gamma\Big(\frac{j}{2}+ t +\frac{1}{2}\Big).
\eeaa
It follows that
\begin{eqnarray}\lbl{volleyball}
\frac{\Gamma_p(\frac{n}{2} + t)}{\Gamma_p(\frac{n}{2})}=\prod_{j=n-p}^{n-1}\frac{\Gamma\big(\frac{j}{2} + t +\frac{1}{2}\big)}{\Gamma\big(\frac{j}{2}+\frac{1}{2}\big)}=\prod_{j=n-p+1}^{n}\frac{\Gamma\big(\frac{j}{2} + t \big)}{\Gamma\big(\frac{j}{2}\big)}.
\end{eqnarray}
This implies
\beaa
\frac{\Gamma_p(\frac{n}{2} + t)}{\Gamma_p(\frac{n}{2})}=\frac{\Gamma\big(\frac{n}{2} + t \big)}{\Gamma\big(\frac{n}{2}\big)}\cdot \frac{\Gamma\big(\frac{n-p}{2} \big)}{\Gamma\big(\frac{n-p}{2} + t\big)}\cdot \prod_{j=n-p}^{n-1}\frac{\Gamma\big(\frac{j}{2} + t \big)}{\Gamma\big(\frac{j}{2}\big)}.
\eeaa
Now, replacing ``$t$" in Proposition \ref{pizza} with ``$-t$" we then obtain
\begin{eqnarray*}
 \log \prod_{j=n-p}^{n-1}\frac{\Gamma(\frac{j}{2}+t)}{\Gamma(\frac{j}{2})}= pt(\log n - 1- \log 2) + r_n^2\Big(t^2 - (p-n+1.5)t\Big) + o(1)
\end{eqnarray*}
as $n\to\infty.$ On the other hand, from Lemma \ref{card},
\begin{eqnarray*}
\log \Big[\frac{\Gamma(\frac{n}{2}+t)}{\Gamma(\frac{n}{2})}\cdot \frac{\Gamma(\frac{n-p}{2})}{\Gamma(\frac{n-p}{2} +t)}\Big] = r_n^2t + o(1)
\end{eqnarray*}
as $n\to\infty.$ Combining the last three equalities, we have
\beaa
\log \frac{\Gamma_p(\frac{n}{2} + t)}{\Gamma_p(\frac{n}{2})}=pt(\log n - 1-\log 2)  + r_n^2\Big(t^2 - (p-n+ \frac{1}{2})t\Big) + o(1)
\eeaa
as $n\to\infty.$  Similarly,
\begin{eqnarray*}
\log \frac{\Gamma_p(\frac{n}{2} + s)}{\Gamma_p(\frac{n}{2})}
 = ps(\log n -1-\log 2) + r_n^2\Big[s^2 -\Big(p-n+\frac{1}{2}\Big)s\Big] + o(1)
\end{eqnarray*}
as $n\to\infty.$ Taking the difference of the above two assertions, we obtain the desired conclusion.\ \ \ \ \ \ \ \ $\blacksquare$

\subsection{Proof of Theorem \ref{yet}}\lbl{Theorem_yet}

\begin{lemma}(Corollary 8.3.6 from Muirhead (1982))\lbl{page62} Assume $n>p.$ Let $V_n$ be as in (\ref{Target}). Then, under $H_0$ in (\ref{have}), we have
\begin{eqnarray*}
E(V_n^h)=p^{ph}\frac{\Gamma(\frac{mp}{2})}{\Gamma(\frac{mp}{2} + ph)}\cdot \frac{\Gamma_p(\frac{m}{2} + h)}{\Gamma_p(\frac{m}{2})}
\end{eqnarray*}
for $h>-\frac{1}{2}$ where $m=n-1.$
\end{lemma}

\noindent\textbf{Proof of Theorem \ref{yet}}.
%For convenience of notation, set $m=n-1.$  Then the theorem is equivalent to that
%\begin{eqnarray}\lbl{skin}
%\frac{\log V_n + p-\big(m-p-\frac{1}{2}\big)r_m^2}{\sigma_n}\ \ \mbox{converges in distribution to}\ \ N(0, 1)
%\end{eqnarray}
%as $n\to\infty,$ where $r_m:=(-\log (1-\frac{p}{m}))^{1/2}$ for $m>p\geq 1.$
Recall that a sequence of random variables $\{Z_n;\, n\geq 1\}$ converges to $Z$ in distribution as $n\to\infty$ if
\bea\lbl{ESPN}
 \lim_{n\to\infty}Ee^{hZ_n}= Ee^{hZ}<\infty
\eea
for all $h\in (-h_0, h_0),$ where $h_0>0$ is a constant. See, e.g., page 408 from Billingsley (1986). Thus, to prove the theorem, it suffices to show that
 there exists $\delta_0>0$ such that
\begin{eqnarray}\lbl{after}
E\exp\Big\{\frac{\log V_n -\mu_n}{\sigma_n}\,s\Big\} \to e^{s^2/2}
\end{eqnarray}
as $n\to\infty$ for all $|s| <\delta_0.$

Set $m=n-1$ and  $r_x:=(-\log (1-\frac{p}{x}))^{1/2}$ for $x>p.$ By the fact that $x+\log (1-x)<0$ for all $x\in (0, 1),$ we know that $\sigma_n^2>0$ for all $n \geq 3,$  and  $\lim_{n\to\infty}\sigma_n^2= -2y-2\log (1-y)>0$ for $y\in (0,1)$, and $\lim_{n\to\infty}\sigma_n^2= +\infty$ for $y=1.$ Therefore,
\begin{eqnarray*}
\delta_0:=\inf\{\sigma_n;\, n\geq 3\}>0
\end{eqnarray*}
Fix $|s|<\frac{\delta_0}{2}.$ Set $t=t_n=\frac{s}{\sigma_n}.$ Then $\{t_n;\, n\geq 3\}$ is bounded and $|t_n|<\frac{1}{2}$ for all $n\geq 3.$  By Lemma \ref{page62},
\begin{eqnarray}\lbl{shoot}
Ee^{t\log V_n}=EV_n^{t}&=& p^{pt}\frac{\Gamma(\frac{mp}{2})}{\Gamma(\frac{mp}{2} + pt)}\cdot \frac{\Gamma_p(\frac{m}{2} + t)}{\Gamma_p(\frac{m}{2})}
\end{eqnarray}
for all $n\geq 3.$ By Lemma \ref{stirling} for the first case and the assumption $p/m\to y \in (0, 1],$
\begin{eqnarray}\lbl{philosophy}
\log \frac{\Gamma(\frac{mp}{2})}{\Gamma(\frac{mp}{2} + pt)} & = & -\log \frac{\Gamma(\frac{mp}{2} + pt)}{\Gamma(\frac{mp}{2})}\nonumber\\
& = & -pt\log \frac{mp}{2} - \frac{p^2t^2- pt}{mp} + O\Big(\frac{1}{m}\Big)\nonumber\\
& = & -pt\log \frac{mp}{2} - \frac{pt^2}{m} + O\Big(\frac{1}{n}\Big)
\end{eqnarray}
as $n\to\infty$. Notice
\beaa
t^2\cdot \left(-\log (1-\frac{p}{m})\right) & = &  \frac{s^2}{\sigma_n^2}\cdot \left(-\log (1-\frac{p}{m})\right)\\
& \to &
\begin{cases}
\frac{s^2}{2}\cdot\frac{\log (1-y)}{y+\log (1-y)}, & \text{if $y\in (0,1)$;}\\
\frac{s^2}{2}, & \text{if $y=1$}
\end{cases}
\eeaa
as $n\to\infty.$ Thus, $t=O(1/r_m)$ as $n\to\infty.$  By Lemma \ref{by},
\beaa
\log \frac{\Gamma_p(\frac{m}{2} + t)}{\Gamma_p(\frac{m}{2})}=pt(\log m-1-\log 2) + r_m^2\left(t^2-(p-m+\frac{1}{2})t\right) + o(1)
\eeaa
as $n\to\infty.$  This together with (\ref{shoot}) and (\ref{philosophy}) gives that
\beaa
\log Ee^{t\log V_n} &= &  pt\log p + \log \frac{\Gamma(\frac{mp}{2})}{\Gamma(\frac{mp}{2} + pt)} + \log \frac{\Gamma_p(\frac{m}{2} + t)}{\Gamma_p(\frac{m}{2})}\\
& = &  pt\log p -pt\log \frac{mp}{2} - \frac{pt^2}{m} + pt(\log m-1-\log 2)  \\
& & \ \ \ \ \ \ \ \ \ \ \  + r_m^2\left(t^2-(p-m+\frac{1}{2})t\right) + o(1)\\
& = & 2\Big(r_m^2 -\frac{p}{m}\Big) \frac{t^2}{2}+\Big[-p+(n-p-\frac{3}{2})r_m^2 \Big] t + o(1)
\eeaa
as $n\to\infty.$ Reviewing the notation $\mu_n$, $\sigma_n$ and $t=t_n=\frac{s}{\sigma_n},$ the above indicates that
\beaa
 \log E\exp\Big\{\frac{\log V_n}{\sigma_n}s\Big\}=\log Ee^{t\log V_n}=\frac{\sigma_n^2t^2}{2} + \mu_nt +o(1)=\frac{s^2}{2} +\frac{\mu_n}{\sigma_n}s + o(1)
\eeaa
as $n\to\infty$ for all $|s|<\frac{\delta_0}{2}.$ This implies (\ref{after}). The proof is completed.\ \ \ \ \ \ \ \ $\blacksquare$

\subsection{Proof of Theorem \ref{energy}}\lbl{Theorem_energy}
\begin{lemma}(Theorem 11.2.3 from Muirhead (1982))\lbl{menu} Let $p$, $n=N-1$ and $W_n$ be as in (\ref{wired5}). Then, under $H_0$ in (\ref{Chicago}),
\begin{eqnarray}\lbl{quiet}
EW_n^t=\frac{\Gamma_p(\frac{n-1}{2}+t)}{\Gamma_p(\frac{n-1}{2})}\cdot \prod_{i=1}^k\frac{\Gamma_{p_i}(\frac{n-1}{2})}{\Gamma_{p_i}(\frac{n-1}{2}+t)}
\end{eqnarray}
for any $t>-1/2,$ where $\Gamma_p(z)$ is as in (\ref{complex}).
\end{lemma}

\noindent\textbf{Proof of Theorem \ref{energy}}. For convenience, set $m=n-1$. Then we need to prove
\begin{eqnarray}\lbl{kiss}
\frac{\log W_n -\tilde{\mu}_m}{\tilde{\sigma}_m}\ \mbox{converges in distribution to}\ N(0,1)
\end{eqnarray}
as $n\to\infty,$ where
\begin{eqnarray*}
 \tilde{\mu}_m=-r_{m}^2\Big(p-m+\frac{1}{2}\Big)+ \sum_{i=1}^k r_{m,i}^2\Big(p_i-m+\frac{1}{2}\Big)\  \mbox{and}\ \  \tilde{\sigma}_m^2=2r_{m}^2-2\sum_{i=1}^k r_{m,i}^2.
\end{eqnarray*}
%$\tilde{\sigma}_m^2=2r_m^2-2\sum_{i=1}^k r_{m,i}^2$, $r_m=(-\log (1-\frac{p}{m}))^{1/2}$, $r_{m,i}=(-\log %(1-\frac{p_i}{m}))^{1/2}$ for $1\leq i \leq k$ and
%\begin{eqnarray*}
%\tilde{\mu}_m=-r_m^2\Big(p-m+\frac{1}{2}\Big)+ \sum_{i=1}^k r_{m,i}^2\Big(p_i-m+\frac{1}{2}\Big).
%\end{eqnarray*}
First, since $m=n-1>p=p_1+\cdots + p_k$  and $\lim_{n\to\infty}p_i/n=y_i$ for each $1\leq i\leq k$, we know
\begin{eqnarray}\lbl{talk}
\frac{p}{m}=\sum_{i=1}^k\frac{p_i}{m}\to \sum_{i=1}^ky_i:=y\in (0, 1]
\end{eqnarray}
as $n\to\infty.$ Second, it is known $\prod_{i=1}^k(1-x_i) > 1-\sum_{i=1}^kx_i$ for all $x_i\in (0,1),\, 1\leq i \leq k,$ see, e.g., p. 60 from Hardy et al. (1988). Taking the logarithm on both sides and then taking $x_i=p_i/m$, we see that
\begin{eqnarray}\lbl{vita}
\frac{1}{2}\tilde{\sigma}_m^2=r_m^2-\sum_{i=1}^k r_{m,i}^2=\sum_{i=1}^k\log (1-\frac{p_i}{m}) - \log (1-\frac{p}{m})>0
\end{eqnarray}
for all $m\geq 2.$ Now, by the assumptions and (\ref{talk}), it is easy to see
\begin{eqnarray*}
\lim_{p\to\infty}\tilde{\sigma}_m^2 =
\begin{cases}
-2\log (1-y)+2\sum_{i=1}^k\log (1-y_i), & \text{if $y<1;$}\\
+\infty, & \text{if $y=1$.}
\end{cases}
\end{eqnarray*}
By the same argument as in the last inequality in (\ref{vita}), we know the limit above is always positive. Reviewing that $p=p_n$ and $m=n-1>p,$ we then have
\begin{eqnarray*}
\delta_0:=\inf\{\tilde{\sigma}_m;\, m\geq 2\}>0.
\end{eqnarray*}
Fix $|s|<\delta_0/2.$ Set $t=t_m=s/\tilde{\sigma}_m.$ Then $\{t_m;\, m\geq 2\}$ is bounded satisfying $|t_m|<1/2$ for all $m\geq 2.$ In particular, as $n\to\infty,$ we have
\bea\lbl{spree}
t=t_m=O\Big(\frac{1}{r_{m,i}}\Big),\ 1\leq i \leq k,
\eea
thanks to that $\lim_{p\to\infty}r_{m,i}= -\log (1-y_i) \in (0, \infty)$ for $1\leq i \leq k.$ On the other hand, notice
\beaa
\sum_{i=1}^kr_{m,i}^2=-\sum_{i=1}^k\log \Big(1-\frac{p_i}{m}\Big) \to -\sum_{i=1}^k\log (1-y_i)
\eeaa
as $n\to\infty.$ It follows from (\ref{talk}) that
\beaa
\lim_{n\to\infty}\frac{r_m^2}{\tilde{\sigma}_m^2}
 =
\begin{cases}
\frac{1}{2}\cdot\frac{\log (1-y)}{\log (1-y)-\sum_{i=1}^k\log (1-y_i)}, & \text{if $y\in (0,1);$}\\
\frac{1}{2}, & \text{if $y=1$.}
\end{cases}
\eeaa
This implies that
\bea\lbl{dust}
t=\frac{s}{\tilde{\sigma}_m}=O\Big(\frac{1}{r_{m}}\Big)
\eea
as $n\to\infty.$ From Lemma \ref{menu},
\begin{eqnarray}\lbl{tilt}
Ee^{t\log W_n}=EW_n^{t}&=& \frac{\Gamma_p(\frac{m}{2}+t)}{\Gamma_p(\frac{m}{2})}\cdot \prod_{i=1}^k\frac{\Gamma_{p_i}(\frac{m}{2})}{\Gamma_{p_i}(\frac{m}{2}+t)}
\end{eqnarray}
since $|t|=|t_m|<1/2.$ By Lemma \ref{by} and (\ref{dust}),
\begin{eqnarray}\lbl{free}
\log \frac{\Gamma_p(\frac{m}{2} + t)}{\Gamma_p(\frac{m}{2})}
 = pt(\log m -1-\log 2) + r_m^2\Big[t^2 -\Big(p-m+\frac{1}{2}\Big)t\Big] + o(1)
\end{eqnarray}
as $n\to\infty.$ Similarly, by Lemma \ref{by} and (\ref{spree}),
\begin{eqnarray}\lbl{defense}
\log \frac{\Gamma_{p_i}(\frac{m}{2} + t)}{\Gamma_{p_i}(\frac{m}{2})}
 = p_it(\log m -1-\log 2) + r_{m,i}^2\Big[t^2 -\Big(p_i-m+\frac{1}{2}\Big)t\Big] + o(1)
\end{eqnarray}
as $n\to\infty$ for  $1\leq i \leq k.$ Therefore, use the identity $p=p_1+\cdots + p_k$ to have
\begin{eqnarray*}
\log \prod_{i=1}^k \frac{\Gamma_{p_i}(\frac{m}{2} + t)}{\Gamma_{p_i}(\frac{m}{2})}
 &= & \sum_{i=1}^k \log \frac{\Gamma_{p_i}(\frac{m}{2} + t)}{\Gamma_{p_i}(\frac{m}{2})}\\
& = & pt(\log m-1-\log 2) +t^2\sum_{i=1}^k r_{m,i}^2 -t\sum_{i=1}^k r_{m,i}^2\Big(p_i-m+\frac{1}{2}\Big)+o(1)
\end{eqnarray*}
as $n\to\infty.$ This together with (\ref{tilt}) and (\ref{free}) gives
\begin{eqnarray*}
\log EW_n^t
&= & t^2\Big(r_m^2-\sum_{i=1}^k r_{m,i}^2\Big) + t\Big[-r_m^2\Big(p-m+\frac{1}{2}\Big)+ \sum_{i=1}^k r_{m,i}^2\Big(p_i-m+\frac{1}{2}\Big)\Big] + o(1)\\
&= & \frac{s^2}{2} + \frac{\tilde{\mu}_m}{\tilde{\sigma}_m}s +o(1)
\end{eqnarray*}
as $n\to\infty$ by the definitions of $\tilde{\mu}_m$ and $\tilde{\sigma}_m$ as well as the fact $t=s/\tilde{\sigma}_m.$ We then arrive at
\begin{eqnarray*}
E\exp\Big\{\frac{\log W_n -\tilde{\mu}_m}{\tilde{\sigma}_m}s\Big\}=e^{-\tilde{\mu}_ms/\tilde{\sigma}_m}\cdot EW_n^t\to e^{s^2/2}
\end{eqnarray*}
as $n\to\infty$ for all $|s|< \delta_0/2.$ This implies (\ref{kiss}) by using the moment generating function method stated in (\ref{ESPN}).\ \ \ \ \ \ \ \ \ \ $\blacksquare$

\subsection{Proof of Theorem \ref{town}}\lbl{Theorem_town}

%{\red In this part we have changed $N_i$ to $n_i$ and $N$ to $n$.}

\noindent Consider
\bea\lbl{ring}
\lambda_n=\frac{\prod_{i=1}^k|\bd{B}_i|^{n_i/2}}{|\bd{A} + \bd{B}|^{n/2}}.
\eea

\begin{lemma}(Corollary 10.8.3 from Muirhead (1982))\lbl{shark} Let $n_i>p$ for $i=1,2\cdots, k.$ Let $\lambda_n$ be as in (\ref{ring}). Then, under $H_0$ in (\ref{multiple}),
\begin{eqnarray*}
E(\lambda_n^t)=\frac{\Gamma_p(\frac{1}{2}(n-1))}{\Gamma_p(\frac{1}{2}n(1+t)-\frac{1}{2})}\cdot \prod_{i=1}^k\frac{\Gamma_p(\frac{1}{2}n_i(1+t)-\frac{1}{2})}{\Gamma_p(\frac{1}{2}(n_i-1))}
\end{eqnarray*}
for all $t>\max_{1\leq i \leq k}\{\frac{p}{n_i}\}-1,$ where $\Gamma_p(z)$ is as in (\ref{complex}).
\end{lemma}
The restriction $t>\max_{1\leq i \leq k}\{\frac{p}{n_i}\}-1$ comes from the restriction in (\ref{complex}).\\

%\begin{prop}\lbl{climb} Suppose $N_i>p$ and $\lim{p/N_i}=y_i\in (0, 1]$ for each $1\leq i \leq k,$ %then
%\begin{eqnarray*}
%\frac{\log \lambda_N - \mu_N}{N\sigma_N}\ \ \mbox{converges in distribution to} \ \ N(0, 1),
%\end{eqnarray*}
% where
%\begin{eqnarray*}
%& & \mu_N=\frac{1}{4}\Big[2\Big(\sum_{i=1}^kN_i\log N_i-N\log N\Big)p
%+  Nr_N^2(2p-2N+3)-\sum_{i=1}^k N_ir_{N_i}^2(2p-2N_i+3)\Big],\\
%& & \sigma_N^2=\frac{1}{2}\Big(\sum_{i=1}^k\frac{N_i^2}{N^2}r_{N_i}^2-r_N^2\Big)>0
%\end{eqnarray*}
%and $r_x=\left(-\log \left(1-\frac{p}{x}\right)\right)^{1/2}$ for $x>p.$
%\end{prop}
\noindent\textbf{Proof of Theorem \ref{town}.} Review (\ref{school}) and (\ref{ring}). Notice
\bea\lbl{depth}
\log \Lambda_n=\log \lambda_n+ \frac{1}{2}pn\log n -\frac{1}{2}\sum_{i=1}^k pn_i\log n_i.
\eea
%Thus, to prove the theorem, it suffices to prove that
%\begin{eqnarray}\lbl{letter}
%\frac{\log \lambda_N - \mu_N'}{N\sigma_N}\ \ \mbox{converges in distribution to} \ \ N(0, 1),
%\end{eqnarray}
% where
%\begin{eqnarray*}
%\mu_N'=\frac{1}{4}\Big[2\Big(\sum_{i=1}^kN_i\log N_i-N\log N\Big)p
%+  Nr_N^2(2p-2N+3)-\sum_{i=1}^k N_ir_{N_i}^2(2p-2N_i+3)\Big].
%\end{eqnarray*}
Evidently,
\begin{eqnarray}\lbl{trio}
\frac{p}{n}=\frac{1}{\sum_{i=1}^k\frac{n_i}{p}}\to \frac{1}{\sum_{i=1}^k\frac{1}{y_i}}:=y \in (0, 1)
\end{eqnarray}
as $p\to \infty.$ As a consequence,
\begin{eqnarray}\lbl{joy}
r_n^2 \to -\log (1-y) \in (0, \infty)
\end{eqnarray}
as $p\to\infty.$

{\it Step 1.} We show $\sigma_n^2>0$ for all $\min_{1\leq i \leq k}n_i\geq p+1$ and $p\geq 1.$ In fact, let $h(x)=-\log (1-x)$ for $x\in [0,1).$ Then $h(x)$ is convex on $[0, 1).$ Take $x_i=\frac{p}{n_i}$, $\lambda_i=\frac{n_i^2}{n^2}$ for $1\leq i \leq k,$ and $x_{k+1}=0$ and $\lambda_{k+1}=1-\frac{\sum_{i=1}^kn_i^2}{n^2}.$ Since $n_i'=n_i-1$ and $r_x$ is decreasing in $x>p,$ by convexity,
\begin{eqnarray}\lbl{brain}
\sum_{i=1}^k\frac{n_i^2}{n^2}r_{n_i'}^2 \geq \sum_{i=1}^k\frac{n_i^2}{n^2}r_{n_i}^2
=\sum_{i=1}^{k+1}\lambda_ih(x_i)>h\Big(\sum_{i=1}^{k+1}\lambda_ix_i\Big) = h\Big(\frac{p}{n}\Big)=r_n^2,
\end{eqnarray}
where ``$>$", instead of ``$\geq$", comes from the fact that $h(x)$ is strictly convex and $x_1\ne x_{k+1}.$ This says that $\sigma_n^2>0$ for all $\min_{1\leq i \leq k}n_i\geq p+1$ and $p\geq 1.$ Second, we claim
\begin{eqnarray}\lbl{small}
2\sigma_n^2\to 2\sigma^2=
\begin{cases}
\log (1-y)-\sum_{i=1}^k\frac{y^2}{y_i^2}\log (1-y_i)>0, & \text{if $\max_{1\leq i\leq k}y_i<1$;}\\
+\infty, & \text{if $\max_{1\leq i\leq k}y_i=1$}
\end{cases}
\end{eqnarray}
as $p\to\infty.$  In fact, for the second case, noticing $y\in (0, 1)$ by (\ref{trio}), the limit is obviously $+\infty$ since $\lim_{x\to 1-}-\log (1-x)=+\infty.$ On the other hand, by (\ref{trio}), $\lim \frac{n_i}{n}=\frac{y}{y_i}$ for all $i=1,2,\cdots, k.$ Thus, the statement  for the case $\max_{1\leq i\leq k}y_i<1$ in (\ref{small}) follows. Moreover, replacing $\frac{p}{n}$ with  $y$, $\frac{p}{n_i}$ with  $y_i$ and $\frac{n_i^2}{n^2}$ with $\frac{y^2}{y_i^2}$ in (\ref{brain}), respectively, we know that the first limit in (\ref{small}) is positive.

{\it Step 2.} In this step we collect some facts that will be used later. Fix $s$ such that $|s|< \frac{\sigma}{2y}.$ Set $t=t_n=\frac{s}{n\sigma_n}.$ We claim that
\begin{eqnarray}
& & t>\max_{1\leq i \leq k}\{\frac{p}{n_i}\}-1\ \ \ \mbox{provided $p$ is sufficiently large};\ \ \ \ \ \ \ \ \ \ \ \ \ \ \lbl{soft}\\
& & \left|\frac{nt-1}{2}\right| \vee \frac{1}{2} =O\Big(\frac{1}{r_n}\Big)\ \ \mbox{and}\ \  \frac{n_it}{2} =O\Big(\frac{1}{r_{n_i-1}}\Big)\lbl{clown}
\end{eqnarray}
as $p\to\infty$ for $i=1,2\cdots, k.$

First, the assumption $\min_{1\leq i \leq k}n_i> p+1$ implies that
\begin{eqnarray}\lbl{night}
\max_{1\leq i \leq k}\{\frac{p}{n_i}\}-1< \frac{p}{p+1}-1=-\frac{1}{p+1}.
\end{eqnarray}
Further, for $|s|< \frac{\sigma}{2y},$ we know $s>-\frac{\sigma}{2y}.$ Moreover, $-\frac{n}{p+1}\sigma_n \to -\frac{\sigma}{y}$ by (\ref{trio}) and (\ref{small}). These imply that, as $p$ is sufficiently large, $-\frac{n}{p+1}\sigma_n <s,$ or $t=\frac{s}{n\sigma_n}> -\frac{1}{p+1}.$ This together with (\ref{night}) concludes (\ref{soft}).

%We now claim that
%\bea\lbl{clown}
%\left|\frac{Nt-1}{2}\right| \vee \frac{1}{2} =O\Big(\frac{1}{r_N}\Big)\ \ \mbox{and}\ \ %\left|\frac{N_it-1}{2}\right| \vee \frac{1}{2} =O\Big(\frac{1}{r_{N_i}}\Big)
%\eea
%as $p\to\infty$ for $i=1,2\cdots, k.$
Second,  by (\ref{joy}),
\bea\lbl{hikers}
\left|\frac{nt-1}{2}\right| \vee \frac{1}{2} =O(1+nt)=O\Big(1+\frac{1}{\sigma_n}\Big)=O\Big(\frac{1}{r_n}\Big)
\eea
as $p\to\infty$ by (\ref{joy}) and (\ref{small}). We obtain the first identity in (\ref{clown}). Moreover, noticing $n_i\leq n$ and $t=s/(n\sigma_n)$ we have
\bea\lbl{paste}
\frac{n_it}{2} =O(nt)=O\Big(\frac{1}{\sigma_n}\Big)
\eea
as $p\to\infty$. By the definition of $\sigma_n^2$, (\ref{joy}) and the fact  that $\lim_{p\to\infty}\frac{n_i}{n}\in (0, \infty)$ again, we know $r_{n_i'}^2=O(\sigma_n^2)$ as $p\to\infty.$ Then $\frac{1}{\sigma_n}=O(\frac{1}{r_{n_i'}})$ as $p\to\infty.$ This joint with (\ref{paste}) gives that
\bea\lbl{checkups}
\frac{n_it}{2} = O\Big(\frac{1}{r_{n_i'}}\Big)
\eea
as $p\to\infty$.
%It is trivial to check that
%\beaa
%\frac{r_{N_i-1}^2}{r_{N_i}^2}=\frac{\log(1-\frac{p}{N_i-1})}{\log(1-\frac{p}{N_i})}=1+ \frac{\log (1-\frac{1}{N_i-p})-\log (1-\frac{1}{N_i})}{\log %(1-\frac{p}{N_i})}=O(1)
%\eeaa
%as $p\to\infty$ since $\log \frac{1}{2} \leq \log (1-\frac{1}{N_i-p}) < 0$ for all $N_i> p+1$ and $\lim_{p\to\infty}\log (1-\frac{p}{N_i}) = \log (1-y_i),$ which is $-\infty$ when $y_i=1.$
This concludes the second identity in (\ref{clown}).

{\it Step 3.} To prove the theorem,  it is enough to prove
\begin{eqnarray}\lbl{plank}
E\exp\Big\{\frac{\log \Lambda_n - \mu_n}{n\sigma_n}s\Big\} \to e^{s^2/2}
\end{eqnarray}
as $p\to \infty$ for all $|s|< \frac{\sigma}{2y}.$

Recalling $t=\frac{s}{n\sigma_n},$ by Lemma \ref{shark} and (\ref{soft}) we have
\begin{eqnarray}\lbl{kick}
Ee^{t\log \lambda_n}=\frac{\Gamma_p(\frac{1}{2}(n-1))}{\Gamma_p(\frac{1}{2}n(1+t)-\frac{1}{2})}\cdot \prod_{i=1}^k\frac{\Gamma_p(\frac{1}{2}n_i(1+t)-\frac{1}{2})}{\Gamma_p(\frac{1}{2}(n_i-1))}.
\end{eqnarray}
Now,  replacing $t$ by $\frac{1}{2}(nt-1)$ and taking $s=-\frac{1}{2}$  in Lemma \ref{by}, by the first assertion in (\ref{clown}) we get
\begin{eqnarray}\lbl{agape}
& & \log \frac{\Gamma_p(\frac{1}{2}(n-1))}{\Gamma_p(\frac{1}{2}n(1+t)-\frac{1}{2})}\nonumber\\
&= & -\log \frac{\Gamma_p(\frac{n}{2} + \frac{nt-1}{2})}{\Gamma_p(\frac{n}{2}-\frac{1}{2})}\nonumber\\
& = & -\frac{n pt}{2}(\log n -1-\log 2) - r_n^2\Big[\frac{(nt-1)^2}{4} -\frac{1}{4} -\Big(p-n+\frac{1}{2}\Big)\frac{nt}{2}\Big] + o(1)\nonumber\\
& = & \frac{n pt}{2}(1+\log 2-\log n) + \frac{r_n^2}{4}\Big(-n^2t^2 + (2p-2n+3)nt\Big) +o(1)
\end{eqnarray}
as $p\to\infty.$ Recall  $n_i'=n_i-1$ for $1\leq i \leq k.$ By Lemma \ref{by} and the second identity in (\ref{clown}),  we get
\begin{eqnarray*}
 \log \frac{\Gamma_p(\frac{1}{2}n_i(1+t)-\frac{1}{2})}{\Gamma_p(\frac{1}{2}(n_i-1))}
& = & \log \frac{\Gamma_p(\frac{1}{2}n_i' + \frac{n_it}{2})}{\Gamma_p(\frac{1}{2}n_i')}\\
% & = & \log \frac{\Gamma_p(\frac{1}{2}(N_i-1))}{\Gamma_p(\frac{1}{2}N_i(1+t)-\frac{1}{2})}\\
% &= & \log \frac{\Gamma_p(\frac{N_i-1}{2} + \frac{N_it}{2})}{\Gamma_p(\frac{N_i-1}{2})}\\
& = & \frac{n_ipt}{2}(\log n_i'-1-\log 2) + \frac{r_{n_i'}^2}{4}\Big(n_i^2t^2 - (2p-2n_i+3)n_it\Big) +o(1)
\end{eqnarray*}
as $p\to\infty.$  Therefore, this, (\ref{kick}) and (\ref{agape}) say that
\begin{eqnarray*}
\log Ee^{t\log \lambda_n}=\frac{t^2}{4}\Big(\sum_{i=1}^kn_i^2r_{n_i'}^2-n^2r_n^2\Big) & + & \frac{t}{4}\Big[2\Big(\sum_{i=1}^kn_i\log n_i'-n\log n\Big)p\\
&+ & nr_n^2(2p-2n+3)-\sum_{i=1}^k n_ir_{n_i'}^2(2p-2n_i+3)\Big]\\
& + & o(1)
\end{eqnarray*}
as $p\to\infty.$  Combining with (\ref{depth}), we obtain that
\begin{eqnarray}
\log Ee^{t\log \Lambda_n}=\frac{t^2}{4}\Big(\sum_{i=1}^kn_i^2r_{n_i'}^2-n^2r_n^2\Big) & + & \frac{t}{4}\Big[2\Big(\sum_{i=1}^kn_i\log n_i'-n_i\log n_i\Big)p\nonumber\\
&+ & nr_n^2(2p-2n+3)-\sum_{i=1}^k n_ir_{n_i'}^2(2p-2n_i+3)\Big]\nonumber\\
& + & o(1)\lbl{Ben}
\end{eqnarray}
as $p\to\infty.$ Observe $n_i\log n_i'-n_i\log n_i=n_i\log(1-\frac{1}{n_i})$ and
\beaa
n_i\log(1-\frac{1}{n_i})  =n_i\Big(-\frac{1}{n_i}-\frac{1}{2}\frac{1}{n_i^2}+O\Big(\frac{1}{n_i^3}\Big)\Big)=-1-\frac{1}{2n_i} + O\Big(\frac{1}{n_i^2}\Big)
\eeaa
as $p\to\infty.$ Thus,
\beaa
p\sum_{i=1}^k(n_i\log n_i'-n_i\log n_i)  =-k p-\frac{1}{2}\sum_{i=1}^ky_i + o(1)
\eeaa
as $p\to\infty.$ Joining this with (\ref{Ben}), we arrive at
\begin{eqnarray*}
\log Ee^{t\log \Lambda_n}=\frac{t^2}{4}\Big(\sum_{i=1}^kn_i^2r_{n_i'}^2-n^2r_n^2\Big) & + & \frac{t}{4}\Big[-2kp-\sum_{i=1}^ky_i +  nr_n^2(2p-2n+3)\nonumber\\ &- & \sum_{i=1}^k n_ir_{n_i'}^2(2p-2n_i+3)\Big] + o(1)
\end{eqnarray*}
as $p\to\infty.$ By the definitions of $\mu_n$ and $\sigma_n,$ the above implies
\begin{eqnarray*}
\log Ee^{t\log \Lambda_n}=\frac{n^2t^2}{2}\sigma_n^2 + \mu_nt + o(1)
\end{eqnarray*}
as $p\to \infty,$ which is equivalent to
\begin{eqnarray*}
\log Ee^{t(\log \Lambda_n-\mu_n)}=\frac{n^2t^2}{2}\sigma_n^2  + o(1)=\frac{s^2}{2} +o(1)
\end{eqnarray*}
as $p\to\infty$ for any $|s|<\frac{\sigma}{2y}$. This leads to (\ref{plank}) since $t=\frac{s}{n\sigma_n}.$\ \ \ \ \ \ \ \ \ $\blacksquare$

\subsection{Proof of Theorem \ref{eve}}\lbl{Theorem_eve}
Let $\Lambda_n^*$ be as in (\ref{stanford3}). Set
\begin{eqnarray}\lbl{boat}
W_n=\frac{\prod_{i=1}^k|\bd{A}_i|^{(n_i-1)/2}}{|\bd{A}|^{(n-k)/2}}=\Lambda_n^*\cdot \frac{\prod_{i=1}^k(n_i-1)^{(n_i-1)p/2}}{(n-k)^{(n-k)p/2}}.
\end{eqnarray}
We have the following result.
\begin{lemma}(p. 302 from Muirhead (1982))\lbl{foot} Assume $n_i>p$ for $1\leq i \leq k.$ Under $H_0$ in (\ref{bang}),
\begin{eqnarray}\lbl{decaf}
E\left[(W_n)^h\right]=\frac{\Gamma_p\left(\frac{1}{2}(n-k)\right)}{\Gamma_p\left(\frac{1}{2}(n-k)(1+h)\right)}\cdot \prod_{i=1}^k\frac{\Gamma_p\left(\frac{1}{2}(n_i-1)(1+h)\right)}{\Gamma_p\left(\frac{1}{2}(n_i-1)\right)}
\end{eqnarray}
for all $h>\max_{1\leq i \leq k}\frac{p-1}{n_i-1}-1,$ where $\Gamma_p(x)$ is defined as in (\ref{complex}).
\end{lemma}
The condition ``$h>\max_{1\leq i \leq r}\frac{p-1}{n_i-1}-1$" is imposed in the above lemma because, by the definition of (\ref{complex}), the following inequalities are needed:
\beaa
%& & \frac{(n_i-1)(1+h)}{2} > \frac{p-1}{2}\ \ \mbox{and}\ \ \  \frac{(n_i-1)}{2} >  \frac{p-1}{2};\ \ \ \ \ \ \ \ \ \ \ \\
& & \frac{(n-k)(1+h)}{2} > \frac{p-1}{2}\ \ \mbox{and} \ \ \ \frac{n-k}{2} >  \frac{p-1}{2};\ \ \ \ \ \ \ \ \ \ \ \\
& & \frac{(n_i-1)(1+h)}{2} > \frac{p-1}{2}\ \ \mbox{and}\ \ \  \frac{(n_i-1)}{2} >  \frac{p-1}{2}
\eeaa
for each $1\leq i \leq k.$ These are obviously satisfied if $n_i>p$ for $1\leq i \leq k$ and $h>\max_{1\leq i \leq k}\frac{p-1}{n_i-1}-1$ (noting that $n=\sum_{i=1}^kn_i>n_i$ for each $i$).\\

\noindent\textbf{Proof of Theorem \ref{eve}}. According to (\ref{boat}), write
\beaa
\log\Lambda_n^*=\log W_n+ \frac{(n-k)p}{2}\log (n-k)- \sum_{i=1}^k\frac{(n_i-1)p}{2}\log (n_i-1).
\eeaa
To prove the theorem, it is enough to show
\begin{eqnarray}\lbl{nun}
\frac{\log W_n - \mu_n'}{(n-k)\sigma_n}\ \ \mbox{converges in distribution to}\ \ N(0, 1)
\end{eqnarray}
as $p\to\infty$, where
\begin{eqnarray}\lbl{BigJack}
 \mu_n' & = & \mu_n +\sum_{i=1}^k \frac{(n_i-1)p}{2} \log ( n_i-1 ) -\frac{(n-k)p}{2} \log (n-k).
\end{eqnarray}
Equation (\ref{nun}) can be proved through the following three steps:

\noindent\underline{\it Step 1}. Let
\begin{eqnarray} \label{y_define}
y = \frac{1}{\sum_{i=1}^k y_i^{-1}}.
\end{eqnarray}
Then, $y\in (0,1)$ and
\begin{eqnarray}\lbl{humility}
\frac{p}{n-k}=\frac{p}{\sum_{i=1}^k (n_i-1)} \to \frac{1}{\sum_{i=1}^k y_i^{-1}} = y.
\end{eqnarray}
 We first show $\sigma_n^2 > 0$. In fact, let $\eta(x)=-\log (1-x)$ for $x\in [0,1).$ Then $\eta(x)$ is strictly convex on $[0, 1).$ Recall that $n=n_1+\cdots + n_k$. Take $x_i=p/(n_i-1)$ and $\lambda_i=(n_i-1)^2/(n-k)^2$ for $1\leq i \leq k$ and $x_{k+1}=0$ and $\lambda_{k+1}=1-\sum_{i=1}^k \lambda_i$. Then, by the strict convexity of $\eta(x)$,
\begin{eqnarray}\lbl{brain1}
-\sum_{i=1}^{k}\Big(\frac{n_i-1}{n-k}\Big)^2\log (1-\frac{p}{n_i-1})
& = & \sum_{i=1}^{k+1}\lambda_i \eta(x_i) >  \eta\Big(\sum_{i=1}^{k+1}\lambda_ix_i\Big) \nonumber \\
& = & \eta\Big( \frac{p}{n-k} \Big) = -\log \Big(1 - \frac{p}{n-k} \Big)
\end{eqnarray}
where the "$>$" holds since $x_{k+1}\ne x_1$. This says that $\sigma_n^2>0$ for all $\min_{1 \leq i \leq k} n_i > 1 + p$ and $p \geq 1$. Secondly, we claim
\begin{eqnarray}\lbl{small1}
2\sigma_n^2\to 2\sigma^2=
\begin{cases}
\log (1-y)-\sum_{i=1}^k\frac{y^2}{y_i^2}\log (1-y_i)>0, & \text{if $\max_{1\leq i\leq k}y_i<1$;}\\
+\infty, & \text{if $\max_{1\leq i\leq k}y_i=1$}
\end{cases}
\end{eqnarray}
as $p\to\infty.$  In fact, the limit in (\ref{small1}) for the case $\max_{1\leq i\leq k}y_i<1$ follows since $n_i/n \to y/y_i$ for all $i=1,2,\cdots, k$ and $-\log (1-x)$ is a continuous function for $x \in (0,1)$. Moreover, replacing $(n_i-1)^2/(n-k)^2$ with $y^2/y_i^2$, $-\log \left[1-p/(n_i-1)\right]$ with  $-\log(1-y_i)$, and $p/(n-k)$ with  $y$ in (\ref{brain1}), respectively, we obtain $\sigma^2> 0$ as $\max_{1\leq i\leq k}y_i<1$. For the second case, we know that one of the $y_i$'s  is equal to 1 and $y\in (0, 1)$. Hence the limit is obviously $+\infty$.

\noindent\underline{\it Step 2}. We will make some preparation for the key part in Step 3. Fix $s$ such that $|s|< \sigma/(2y).$ Set $t=t_n=s/[(n-k) \sigma_n]$. We claim that
\begin{eqnarray}
& & t>\max_{1\leq i \leq k}\Big\{\frac{p-1}{n_i-1}\Big\}-1\ \ \mbox{as $p$ is sufficiently large},\ \ \ \ \ \ \ \ \ \ \ \ \ \lbl{soft1}\\
& & \frac{(n-k)t}{2}=O\Big(\frac{1}{r_{n-k}}\Big)\ \ \mbox{and}\ \ \frac{(n_i-1)t}{2}=O\Big(\frac{1}{r_{n_i-1}}\Big),\lbl{sweater}
\end{eqnarray}
as $p\to\infty$ for $1\leq i \leq k$, where $r_{n-k}=\{-\log (1-\frac{p}{n-k})\}^{1/2}$ and $r_{n_i-1}=\{-\log (1-\frac{p}{n_i-1})\}^{1/2}$.  First, the assumption $\min_{1\leq i \leq k}n_i > 1 + p$ implies that
\begin{eqnarray}\lbl{night1}
\max_{1\leq i \leq k}\Big\{\frac{p-1}{n_i-1}\Big\}-1< \frac{p-1}{p}-1=-\frac{1}{p}.
\end{eqnarray}
Further, since $|s|< \sigma/(2y)$, we see $s>-\sigma/(2y)$. Moreover, $-(n-k)\sigma_n/p \to -\sigma/y$ as $p\to\infty$ by (\ref{humility}) and (\ref{small1}). These imply that, as $p$ is sufficiently large, $-(n-k)\sigma_n/p <s$, or $t=s/[(n-k) \sigma_n]> -1/p$. This together with (\ref{night1}) concludes (\ref{soft1}).

Secondly, since $t=t_n=s/[(n-k) \sigma_n]$, we know from (\ref{small1}) that $\{(n-k)t_n\}$ is bounded. Then the first assertion in (\ref{sweater}) follows since $r_{n-k}^2\to -\log (1-y)\in (0, \infty)$ as $p\to\infty$. Now, fix an $i \in \{1, \cdots, k\}$. Easily, by (\ref{small1}),
\bea\lbl{dense}
\frac{r_{n_i-1}^2}{\sigma_n^2} \to
\begin{cases}
\frac{-\log (1-y_i)}{\sigma^2}, & \text{if $\max_{1\leq j\leq k}y_j<1$,}\\
0, & \text{if $\max_{1\leq j\leq k}y_j=1$ and $y_i<1$}
\end{cases}
\eea
as $p\to\infty.$ Now assume $y_i=1$ for some $1\leq i \leq k.$ By the definition of $\sigma_n^2,$ we see that
\beaa
2\sigma_n^2\geq -r_{n-k}^2 + \left(\frac{n_i-1}{n-k}\right)^2 r_{n_i-1}^2\to +\infty
\eeaa
 as $p\to\infty.$ Therefore, use the facts that $r_{n-k}^2\to -\log (1-y)\in (0, \infty)$ and $(n_i-1)/(n-k)\to y/y_i\in (0, \infty)$  to have
\beaa
\frac{(n_i-1)t}{2}=\frac{s}{2}\cdot \frac{n_i-1}{n-k}\cdot \frac{1}{\sigma_n} =O\Big(\frac{1}{r_{n_i-1}}\Big)
\eeaa
as $p\to\infty$. Combining this with (\ref{dense}) we see that
\begin{eqnarray*}
\frac{(n_i-1)t}{2}=O\Big(\frac{1}{r_{n_i-1}}\Big)
\end{eqnarray*}
as $p\to\infty$ for any $y_i\in (0, 1]$. This gives the second assertion in (\ref{sweater}).

\noindent\underline{\it Step 3}.
To prove the theorem, from (\ref{ESPN}) and (\ref{nun}) it suffices to prove
\begin{eqnarray}\lbl{plank1}
E\exp\Big\{\frac{\log W_n - \mu_n'}{(n-k)\sigma_n}s\Big\} \to e^{s^2/2}
\end{eqnarray}
as $p\to \infty$ for all $|s|< \sigma/(2y)$. Recall $t=t_n=s/[(n-k) \sigma_n]$. By Lemma \ref{foot} and (\ref{soft1}),
%Then according to the moment result in (\ref{multiple_covariance_moments}), apparently,
\begin{eqnarray}
\log E \exp \Big\{\frac{\log W_n}{(n-k)\sigma_n}s\Big\} & = & \log E\left[W_n^t\right] \nonumber \\
& = & \log \frac{\Gamma_p \big[ \frac{1}{2}(n-k) \big]} {\Gamma_p \left[ \frac{1}{2}(n-k)(1+t)\right] } +
\sum_{i=1}^k \log \frac{\Gamma_p \left[ \frac{1}{2}(n_i-1)(1+t)\right]} {\Gamma_p \left[ \frac{1}{2}(n_i-1) \right]} \nonumber
\end{eqnarray}
as $p$ is sufficiently large. Using Lemma \ref{by} and the first assertion of (\ref{sweater}), we obtain
\begin{eqnarray}\lbl{vase}
 \log \frac{\Gamma_p \left[ \frac{1}{2}(n-k) \right]} {\Gamma_p \left[ \frac{1}{2}(n-k)(1+t)\right] }
& = & \log \frac{\Gamma_p \left[ \frac{1}{2}(n-k) \right]} {\Gamma_p \left[ \frac{1}{2}(n-k)+ \frac{1}{2}(n-k) t \right]} \nonumber \\
& = & -\frac{(n-k)pt}{2} \left[ \log (n-k) - 1 - \log 2 \right] \nonumber \\
& & - r_{n-k}^2 \Big[ \frac{(n-k)^2t^2}{4} - \frac{(p-n+k+0.5)(n-k)t}{2} \Big] + o(1) \nonumber
\end{eqnarray}
as $p\to\infty$. Similarly, by the Lemma \ref{by} and the second assertion of (\ref{sweater}), we have
 \begin{eqnarray*}
\log \frac{\Gamma_p \big[ \frac{1}{2}(n_i-1)(1+t)\big]} {\Gamma_p \big[ \frac{1}{2}(n_i-1) \big]}
 & = & \frac{(n_i-1)pt}{2}\left[\log (n_i-1) -1-\log 2\right] \\
 & & \ \ \ \ + r_{n_i-1}^2\Big[\frac{(n_i-1)^2t^2}{4} -\frac{(p-n_i+1.5)(n_i-1)t}{2}\Big] +o(1)
\end{eqnarray*}
as $p\to\infty$. Take sum over all $i$ to have
 \begin{eqnarray*}
& & \sum_{i=1}^k \log \frac{\Gamma_p \big[ \frac{1}{2}(n_i-1)(1+t)\big]} {\Gamma_p \big[ \frac{1}{2}(n_i-1) \big]} \\
& = & \frac{pt}{2}\Big[\sum_{i=1}^k (n_i-1)\log (n_i-1)\Big]-\frac{(n-k)pt}{2}(1+\log 2) +\frac{t^2}{4}\Big[\sum_{i=1}^{k}(n_i-1)^2r_{n_i-1}^2\Big] \\
& & \ -\frac{t}{2}\sum_{i=1}^k \left(p-n_i+1.5\right)(n_i-1)r_{n_i-1}^2 + o(1)
\end{eqnarray*}
as $p\to\infty$. Therefore,
\begin{eqnarray*}
\log Ee^{t\log W_n}
& = & \frac{t^2}{4} \Big[\sum_{i=1}^k (n_i-1)^2 r_{n_i-1}^2 - (n-k)^2 r_{n-k}^2 \Big] \\
& & + \frac{t}{2} \Big[ (p-n+k+0.5)(n-k)r_{n-k}^2 - \sum_{i=1}^k (p-n_i+1.5)(n_i-1) r_{n_i-1}^2 \Big] \\
& & + \frac{pt}{2} \Big[ \sum_{i=1}^k(n_i-1)\log (n_i-1) - (n-k) \log (n-k) \Big] + o(1) \\
& = & \frac{t^2}{2} (n-k)^2 \sigma_n^2 + \mu_{n}^\prime t + o(1)
\end{eqnarray*}
as $p\to\infty$ where $\mu_n'$ is as in (\ref{BigJack}). Since $t=t_n=s/[(n-k) \sigma_n]$, we know
\begin{eqnarray*}
\log E\exp\Big\{\frac{\log W_n - \mu_n'}{(n-k)\sigma_n}s\Big\} = \log Ee^{t\log W_n} -\mu_n't\to \frac{s^2}{2}
\end{eqnarray*}
as $p\to\infty$ for all $|s| < \sigma/(2y).$ This leads to (\ref{plank1}). \ \ $\blacksquare$\\

\subsection{Proof of Theorem \ref{Ruth}}\lbl{Theorem_Ruth}

%The following is from Theorem 8.5.1 and Corollary 8.5.4 in \cite{Muirhead1982}.
\begin{lemma}(Theorems 8.5.1 and 8.5.2 and Corollary 8.5.4 from Muirhead (1982))\lbl{mall} Assume $n>p.$ Then the LRT statistic for testing  $H_0$ in (\ref{hypothesis_mean_covariance_specified_value_b}) is given by
\beaa
\Lambda_n=\Big(\frac{e}{n}\Big)^{np/2}|\bd{A}|^{n/2}\cdot \exp\Big\{-\frac{1}{2}tr(\bd{A})-\frac{1}{2}n\overline{\bd{x}}'\overline{\bd{x}})\Big\}
\eeaa
is unbiased, where $\bd{A}$ is as in (\ref{oar}). Further, assuming $H_0$ in (\ref{hypothesis_mean_covariance_specified_value_b}), we have
\beaa
E(\Lambda_n^t)=\Big(\frac{2e}{n}\Big)^{npt/2}(1+t)^{-np(1+t)/2}\frac{\Gamma_p(\frac{n(1+t)-1}{2})}{\Gamma_p(\frac{n-1}{2})}
\eeaa
for any $t>\frac{p}{n}-1.$
\end{lemma}
The range ``$t>\frac{p}{n}-1$"  follows from the definition of $\Gamma_p(z)$ in (\ref{complex}).\\

\noindent\textbf{Proof of Theorem \ref{Ruth}}. First, since $\log (1-x)< -x$ for all $x<1,$ we know $\sigma_n^2>0$ for all $n\geq 3.$ Now, by assumption, it is easy to see that
\begin{eqnarray}\lbl{screen}
\lim_{n\to\infty}\sigma_n^2 =
\begin{cases}
-\frac{1}{2}\big[y+\log (1-y)\big], & \text{if $y<1;$}\\
+\infty, & \text{if $y=1$.}
\end{cases}
\end{eqnarray}
Easily, the limit is always positive. Hence,
\begin{eqnarray*}
\delta_0:=\inf\{\sigma_n;\, n\geq 3\}>0.
\end{eqnarray*}
Fix a number $h$ with $|h|< \delta_0,$ then $h>-\delta_0 \geq -\sigma_n$ for all $n\geq 3.$ It follows that
\beaa
\frac{p}{n}-1\leq \frac{n-1}{n}-1=-\frac{1}{n}< \frac{h}{n\sigma_n}
\eeaa
for all $n\geq 3.$   Set $t=t_n=\frac{h}{n\sigma_n}.$ Then the above says that
\begin{eqnarray}\lbl{bulb}
t>\frac{p}{n}-1
\end{eqnarray}
for all $n\geq 3.$ From  (\ref{screen}) we know that $\{t_n;\, n\geq 3\}$ is bounded. By Lemma \ref{mall} and (\ref{bulb}),
\begin{eqnarray}\lbl{hair}
Ee^{t\log \Lambda_n}=E(\Lambda_n^{t})&=& \Big(\frac{2e}{n}\Big)^{npt/2}(1+t)^{-np(1+t)/2}\frac{\Gamma_p(\frac{n(1+t)-1}{2})}{\Gamma_p(\frac{n-1}{2})}
\end{eqnarray}
for all $n\geq 3.$ To prove the theorem, we only need to show
\begin{eqnarray}\lbl{Albert}
E\exp\Big\{\frac{\log \Lambda_n - \mu_n}{n\sigma_n}h\Big\} \to e^{h^2/2}
\end{eqnarray}
as $n\to\infty$ for all $h$ with $|h|< \delta_0.$ Let $r_n=(-\log (1-\frac{p}{n}))^{1/2}$ for all $n\geq 3.$ From the definition of $\sigma_n^2$ and (\ref{screen}), it is evident that
\bea\lbl{silver}
%\Big(\Big|\frac{nt-1}{2}\Big|\vee \frac{1}{2}\Big)^2\cdot r_n^2=O\left((n|t| +1)^2 r_n^2\right)=O\left(\frac{r_n^2}{\sigma_n^2}\right)=O(1)
\frac{nt}{2}=\frac{h}{2\sigma_n}=O\Big(\frac{1}{r_{n-1}}\Big)
\eea
as $n\to\infty.$   Set $m=n-1$. Take $s=0$ and replace $t$ with $\frac{nt}{2}$  in Lemma \ref{by}, we obtain from (\ref{silver}) that
\begin{eqnarray}\lbl{flower}
\log \frac{\Gamma_p(\frac{n(1+t)-1}{2})}{\Gamma_p(\frac{n-1}{2})}
%\nonumber\\
& = & \log \frac{\Gamma_p(\frac{m}{2}+\frac{nt}{2})}{\Gamma_p(\frac{m}{2})}\nonumber\\
%& = & \log \frac{\Gamma_p(\frac{n}{2} + \frac{nt-1}{2})}{\Gamma_p(\frac{n}{2}-\frac{1}{2})}\nonumber\\
& = & \frac{npt}{2}(\log m -1-\log 2) + r_m^2\Big[\frac{n^2t^2}{4}  -\frac{1}{2}\Big(p-m+\frac{1}{2}\Big)nt\Big] +o(1)\nonumber\\
&= & \frac{npt}{2}(\log m -1-\log 2) + \frac{n^2r_m^2}{4}t^2 +\frac{1}{4}\Big(2n-2p-3\Big)nr_m^2t + o(1) \nonumber\\
& &
\end{eqnarray}
as $n\to\infty.$ Use the fact $\log (1+s)=s-\frac{s^2}{2} + o(s^3)$ as $s\to 0$ to have
\beaa
(1+t)\log (1+t)=(1+t)\Big(t-\frac{t^2}{2} + o(t^3)\Big)=t+\frac{t^2}{2} + O(t^3)
\eeaa
as $n\to\infty$ since $\lim_{n\to\infty}t=\lim_{n\to\infty}t_n= 0.$ It follows that
\beaa
\log (1+t)^{-np(1+t)/2}=-\frac{1}{2}np(1+t)\log (1+t)=-\frac{1}{2}np\Big(t+\frac{t^2}{2}\Big) + O(npt^3)
\eeaa
as $n\to\infty.$ Thus, by (\ref{hair}),
\bea
& & \log E(\Lambda_n^t)\nonumber\\
 &= & \frac{npt}{2}\log \frac{2e}{n}-\frac{1}{2}np\Big(t+\frac{t^2}{2}\Big) + \frac{npt}{2}(\log m -1-\log 2)\lbl{june} \\
 & & \ \ \ \ \ \ \ \ \ \ \ \ \ \ \ + \frac{n^2r_m^2}{4}t^2 +\frac{1}{4}\Big(2n-2p-3\Big)nr_m^2t+ O(npt^3) +o(1)\nonumber\\
 &= & \frac{1}{4}(n^2r_m^2-np)t^2 + \frac{npt}{2}\log (1- \frac{1}{n}) +\frac{1}{4}\Big[(2n-2p-3)nr_m^2 -2np\Big]t +o(1)\ \ \ \ \ \ \lbl{fiber}
\eea
where the sum of the first and third terms in (\ref{june}) gives the second term in (\ref{fiber}), and $O(npt^3)=o(1)$ as $n\to\infty$ by the definition of $t$ and (\ref{screen}). Now,
\beaa
% & & \frac{npt}{2}\log (1- \frac{1}{n})=\frac{npt}{2}\Big(-\frac{1}{n} + O\Big(\frac{1}{n^2}\Big)\Big)=-\frac{pt}{2} + o(1)\ \mbox{and}\\
 \frac{1}{4}(n^2r_m^2-np)t^2 & = & \frac{1}{4}n^2\Big(r_m^2-\frac{p}{m}\Big)t^2 + \frac{1}{4}n^2\Big(\frac{p}{m}-\frac{p}{n}\Big)t^2\\
& = & \frac{1}{4}n^2\Big(r_m^2-\frac{p}{m}\Big)t^2 + o(1)
\eeaa
since $\frac{1}{4}n^2(\frac{p}{m}-\frac{p}{n})t^2=O(pt^2)=o(1)$ as $n\to\infty$  by the definition of $t$ and (\ref{screen}). Also,
\beaa
 \frac{npt}{2}\log (1- \frac{1}{n})=\frac{npt}{2}\Big(-\frac{1}{n} + O\Big(\frac{1}{n^2}\Big)\Big)=-\frac{pt}{2} + o(1)
\eeaa
as $n\to\infty.$  Joining the above two assertions and (\ref{fiber}), recalling the definitions of $\mu_n$ and $\sigma_n$, we get
\beaa
\log E(\Lambda_n^t) & = & \frac{1}{2}n^2t^2\cdot \frac{1}{2}\Big(r_m^2-\frac{p}{m}\Big)+\frac{1}{4}\Big[(2n-2p-3)nr_m^2 -2(n+1)p\Big]t +o(1)\\
& = & \frac{n^2\sigma_n^2t^2}{2} + \mu_nt + o(1)=\frac{h^2}{2} + \mu_nt + o(1)
\eeaa
as $n\to\infty$ for all $h$ with $|h|<\delta_0$ since $t=t_n=\frac{h}{n\sigma_n}.$ Therefore, we eventually conclude
\begin{eqnarray*}
\log E\exp\Big\{\frac{\log \Lambda_n - \mu_n}{n\sigma_n}h\Big\}=\log E(\Lambda_n^t)-\mu_nt \to \frac{h^2}{2}
\end{eqnarray*}
as $n\to\infty$ for all $|h|<\delta_0,$ which is (\ref{Albert}). The proof is completed.\ \ \ \ \ \ \ $\blacksquare$

\subsection{Proof of Theorem \ref{season}}\lbl{Theorem_season}

%According to (9) on p. 150 from \cite{Muirhead1982} or (48) on p. 492 from \cite{Wilks}, for all $n>p\geq 1$,
\begin{lemma}\lbl{cookup} Let $\bd{R}_n$ be the correlation matrix with the density function as in (\ref{adopt}). Assume $n-5\geq p\geq 2.$ Then,
\begin{eqnarray}\lbl{go}
E[|\bd{R}_n|^t]=\Big[\frac{\Gamma(\frac{n-1}{2})}{\Gamma(\frac{n-1}{2}+t)}\Big]^p\cdot \frac{\Gamma_p(\frac{n-1}{2}+t)}{\Gamma_p(\frac{n-1}{2})}
\end{eqnarray}
for all $t\geq -1.$
\end{lemma}
In the literature such as Muirhead (1982)  and Wilks ({1932}), the above formula is only valid for integer $t\geq -1.$ The above lemma says that it is actually  true for all real number $t\geq -1.$\\

\noindent\textbf{Proof of Lemma \ref{cookup}}.
%Since $\Gamma_1(z)=\Gamma(z)$ for all $z\in \mathbb{C}$ and $\bd{R}_n=(1)_{1\times 1}$ for $p=1,$  we know that (\ref{go}) is obviously true %for $p=1.$ So, without loss of generality, assume $p\geq 2.$
Recall (\ref{playful}),  $\bd{R}_n$ is a $p\times p$ non-negative definite matrix and each of its entries takes value in $[-1, 1],$ thus the determinant $|\bd{R}_n|\leq p!.$ By (9) on p. 150 from Muirhead (1982)  or (48) on p. 492 from Wilks ({1932}),
\bea\lbl{youngest}
E\big[|\bd{R}_n|^k\big]= \Big[\frac{\Gamma(\frac{n-1}{2})}{\Gamma(\frac{n-1}{2}+k)}\Big]^p\cdot \frac{\Gamma_p(\frac{1}{2}(n-1)+k)}{\Gamma_p(\frac{1}{2}(n-1))}
\eea
for any integer $k$ such that $\frac{n-1}{2}+k>\frac{p-1}{2}$ by (\ref{complex}), which is equivalent to that $n-p>-2k.$ By the given condition, $n-p\geq 5.$ Thus, (\ref{youngest}) holds for all $k\geq -2.$ In particular,  $E\big[|\bd{R}_n|^{-2}\big]< \infty.$
%By the inequality $|\log x| \leq x+ x^{-1}$ for all $x>0$, we have $E|\log |\bd{R}_n||< \infty.$
Since $|\bd{R}_n|$ is bounded, this implies that
\bea\lbl{haoma}
E\big[|\bd{R}_n|^{-1}\big]<\infty\ \mbox{and}\ E\big[|\bd{R}_n|^{-1}|\log |\bd{R}_n|\,|\big]<\infty.
\eea
Now set $Z=-\log (|\bd{R}_n|/p!).$ Then $P(Z>0)=1$, $E(Ze^{Z})<\infty$ by (\ref{haoma}) and
\bea\lbl{rememberme}
h_1(z):=(p!)^{-(z-1)}\cdot E\big[|\bd{R}_n|^{z-1}\big]= Ee^{-(z-1)Z}
\eea
for all $\mbox{Re}(z)\geq 0.$ It is not difficult to check that $\frac{d}{dz}\,(Ee^{-(z-1)Z})=-E\big[Ze^{-(z-1)Z}\big]$ for all $\mbox{Re}(z)\geq 0.$  Further, by (\ref{haoma}) again, $|h_1(z)|\leq p!\cdot E\big[|\bd{R}_n|^{-1}|\big]<\infty$ on $\{z\in\mathbb{C};\,\mbox{Re}(z)\geq 0\}.$ Therefore, $h_1(z)$ is analytic and bounded on $\{z\in\mathbb{C};\,\mbox{Re}(z)\geq 0\}.$ Define
\bea\lbl{partyboy}
h_2(z)=(p!)^{-(z-1)}\cdot\Big[\frac{\Gamma(\frac{n-1}{2})}{\Gamma(\frac{n-1}{2}+z-1)}\Big]^p\cdot \frac{\Gamma_p(\frac{n-1}{2}+z-1)}{\Gamma_p(\frac{1}{2}(n-1))}
\eea
for $\mbox{Re}(z)\geq 0.$ By the Carlson uniqueness theorem (see, for example Theorem 2.8.1 on p. 110 from Andrews et al. (1999)), if we know that $h_2(z)$ is also bounded and analytic on $\{z\in\mathbb{C};\,\mbox{Re}(z)\geq 0\},$ since $h_1(z)=h_2(z)$ for all $z=0,1,2,\cdots$, we obtain that $h_1(z)=h_2(z)$ for all $\mbox{Re}(z)\geq 0.$ This implies our desired conclusion. Thus,  we only need to check that $h_2(z)$  is bounded and analytic on $\{z\in\mathbb{C};\,\mbox{Re}(z)\geq 0\}.$ To do so, review (\ref{complex}), it suffices to show
\bea
h_3(z):=\prod_{i=2}^p\frac{\Gamma(\frac{n-i}{2}+z-1)}{\Gamma(\frac{n-1}{2}+z-1)}
\eea
is bounded and analytic on $\{z\in\mathbb{C};\,\mbox{Re}(z)\geq 0\}.$ Noticing $\frac{3}{2}\leq \frac{n-i}{2}-1\leq \frac{n-2}{2}-1$ for all $2\leq i \leq p$, to show that, it is enough to prove
\bea\lbl{babyboy}
h(z):=\frac{\Gamma(\alpha +z)}{\Gamma(\beta+z)}
\eea
is bounded and analytic on $\{z\in \mathbb{C};\, \mbox{Re}(z)\geq 0\}$ for all fixed $\beta>\alpha>0.$ This is confirmed by Lemma \ref{latex} in the Appendix.\ \ \ \ \ \ $\blacksquare$\\

\noindent\textbf{Proof of Theorem \ref{season}}. First, since $\log (1-x)< -x$ for all $x<1,$ we know $\sigma_n^2>0$ for all $n\geq 3$ since $n-1>p\geq 1$ for all $n\geq 3$ by the assumption. Now, from the given condition, it is easy to see
\begin{eqnarray}\lbl{chase}
\lim_{n\to\infty}\sigma_n^2 =
\begin{cases}
-2\big[y+\log (1-y)\big], & \text{if $y<1;$}\\
+\infty, & \text{if $y=1.$}
\end{cases}
\end{eqnarray}
Trivially, the limit is always positive for $y\in (0,1]$. Consequently,
\begin{eqnarray*}
\delta_0:=\inf\{\sigma_n;\, n\geq 3\}>0.
\end{eqnarray*}
To finish the proof, by (\ref{ESPN}) it is enough to show that
\begin{eqnarray}\lbl{nose}
E\exp\Big\{\frac{\log |\bd{R}_n|-\mu_n}{\sigma_n}s\Big\}\to e^{s^2/2}
\end{eqnarray}
as $n\to\infty$  for all $s$ such that $|s|<\delta_0/2.$

\noindent Fix $s$ such that $|s|<\delta_0/2.$ Set $t=s/\sigma_n.$ Then $|t|<1/2$ for all $n\geq 3.$ Thus, by Lemma \ref{stirling} for the second case,
\begin{eqnarray}
\log \Big[\frac{\Gamma(\frac{n-1}{2})}{\Gamma(\frac{n-1}{2}+t)}\Big]^p
&= & -p\cdot\log\frac{\Gamma(\frac{n-1}{2}+t)}{\Gamma(\frac{n-1}{2})}\nonumber\\
&= & -p\left(t\log \frac{n-1}{2}+\frac{t^2-t}{n-1} +O(\frac{1}{n^2})\right)\nonumber\\
& = & -pt[\log (n-1)-\log 2] -\frac{p}{n-1}(t^2-t) + o(1)\lbl{meatball}
%& = & -pt\log n +pt\log 2 -yt^2
\end{eqnarray}
as $n\to\infty.$  Second, it is ease to see that
\beaa
t^2r_{n-1}^2=\frac{r_{n-1}^2}{\sigma_n^2}s^2\to
\begin{cases}
\frac{s^2}{2}\frac{\log (1-y)}{y+\log (1-y)}, & \text{if $y\in (0, 1);$}\\
\frac{s^2}{2}, & \text{if $y=1$}
\end{cases}
\eeaa
as $n\to\infty,$ where  $r_n=(-\log (1-\frac{p}{n}))^{1/2}$ for all $n\geq 3.$ In particular, $t=O(1/r_{n-1})$ as $n\to\infty.$ Therefore, by Lemma \ref{by},
\begin{eqnarray*}
\log \frac{\Gamma_p(\frac{n-1}{2}+t)}{\Gamma_p(\frac{n-1}{2})}=pt(\log (n-1)-1-\log 2) + r_{n-1}^2\left(t^2-(p-n+\frac{3}{2})t\right) + o(1)
\end{eqnarray*}
as $n\to\infty.$ By the given condition, $\bd{R}_n$ has the density function as in (\ref{adopt}). Therefore, from  (\ref{go}) and (\ref{meatball}) we conclude that
\begin{eqnarray*}
\log E[|\bd{R}_n|^t] & = & \Big(r_{n-1}^2-\frac{p}{n-1}\Big)t^2 -\Big[p-\frac{p}{n-1}+(p-n+\frac{3}{2})r_{n-1}^2\Big]t+o(1)\\
& = & \frac{s^2}{2} + \mu_nt +o(1)
\end{eqnarray*}
as $n\to\infty$ since $t=s/\sigma_n$ and $\sigma_n^2=2(r_{n-1}^2-\frac{p}{n-1}).$ This implies that
\begin{eqnarray*}
E\exp\Big\{\frac{\log |\bd{R}_n|-\mu_n}{\sigma_n}s\Big\}=e^{-\mu_nt} E[|\bd{R}_n|^t]\to e^{s^2/2}
\end{eqnarray*}
as $n\to\infty$ for any $|s|< \delta_0/2.$ We get (\ref{nose}).\ \ \ \ \ \ \ \ $\blacksquare$

\subsection{Appendix}\lbl{appendix}
%In this part, by using Euler's formula on the Gamma function $\Gamma(z)$, we give two lemmas needed to prove Lemmas
In this section we give a lemma on the complex analysis needed to prove Lemma \ref{cookup}.

\begin{lemma}\lbl{latex} Let $\Gamma(z)$ be the Gamma function defined on the complex plane $\mathbb{C}.$ Let $\beta>\alpha>0$ be two constants. Then
\beaa
h(z):=\frac{\Gamma(\alpha +z)}{\Gamma(\beta+z)}
\eeaa
is bounded and analytic on $\{z\in \mathbb{C};\, \mbox{Re}(z) \geq 0\}.$
\end{lemma}
\textbf{Proof}. It is known that $\Gamma(z)$ is a meromorphic function and all of its poles are simple poles at $z=0,-1,-2,\cdots$ Also $\Gamma(z)$ has no zeros on the complex plane $\mathbb{C}$ (see, e.g., p. 199 from Ahlfors (1979) or  p. 364 from Gamelin (2001)). Thus, $h(z)$ is analytic for all $\mbox{Re}(z)\geq 0.$ On the other hand,  by Euler's formula (see, e.g., p. 199 from Ahlfors (1979) or p. 363 from Gamelin (2001)),
\bea\lbl{terrific}
\frac{1}{\Gamma(z)}=ze^{\gamma z}\prod_{k=1}^\infty\big(1+\frac{z}{k}\big)e^{-z/k}
\eea
for all $z\in \mathbb{C},$ where $\gamma=0.5772\cdots$ is the Euler constant. Hence,
\bea\lbl{chin}
%\frac{\Gamma(\alpha +z)}{\Gamma(\beta+z)}
h(z)=\frac{z+\beta}{z+\alpha}\cdot e^{\gamma(\beta -\alpha)}\cdot \prod_{k=1}^{\infty}\frac{k+z+\beta}{k+z+\alpha}\cdot e^{-(\beta-\alpha)/k}
\eea
for all $\mbox{Re}(z)\geq 0.$  Since $|k+z+\alpha|\geq k+\alpha$ for all $\mbox{Re}(z)\geq 0,$ we have
\beaa
\Big|\frac{k+z+\beta}{k+z+\alpha}\Big|=\Big|1+\frac{\beta-\alpha}{k+z+\alpha}\Big| \leq 1+\frac{\beta-\alpha}{k+\alpha} \leq \exp\Big\{\frac{\beta-\alpha}{k+\alpha}\Big\}
\eeaa
for all $\mbox{Re}(z)\geq 0.$ Consequently,
\beaa
\Big|\prod_{k=1}^{\infty}\frac{k+z+\beta}{k+z+\alpha}\cdot e^{-(\beta-\alpha)/k}\Big| \leq \exp\Big\{-\sum_{k=1}^{\infty}\frac{(\beta-\alpha)\alpha}{k(k+\alpha)}\Big\} \leq 1
%\leq e^{-\pi^2(\beta-\alpha)\alpha/6}
\eeaa
for all $\mbox{Re}(z)\geq 0$ since $\beta>\alpha>0$.
%by the fact $\sum_{k=1}^{\infty}\frac{1}{k(k+\alpha)} < \infty$
Obviously, $\frac{\alpha +z}{\beta+z}$ is bounded on $\{z\in \mathbb{C};\, \mbox{Re}(z)\geq 0\}.$ By (\ref{chin}) and the first paragraph, we know $h(z)$ is bounded and analytic on $\{z\in \mathbb{C};\, \mbox{Re}(z)\geq 0\}.$\ \ \ \ \ \ $\blacksquare$\\

%\begin{lemma}\lbl{areas} Let $\Gamma(x)$ be the Gamma function defined on $(0, \infty).$ Let $d>c>0$ be two constants. Then, as %$\epsilon\to 0,$
%\beaa
%\sup_{c\leq x \leq d}\Big|\log \frac{\Gamma(x+\epsilon)}{\Gamma(x)}\Big| \to 0.
%\eeaa
%\end{lemma}
%\textbf{Proof}. By (\ref{terrific}),
%\beaa\lbl{Thomson}
%\frac{\Gamma(x+\epsilon)}{\Gamma(x)}=\frac{x}{x+\epsilon}\cdot e^{-\gamma\epsilon}\cdot %\prod_{k=1}^{\infty}\Big(\frac{k+x}{k+x+\epsilon}\cdot e^{\epsilon/k}\Big)
%& \leq & e^{(2c^{-1}+\gamma)|\epsilon|}\cdot \prod_{k=1}^{\infty}\Big(e^{-\epsilon/(k+x+\epsilon)}\cdot e^{\epsilon/k}\Big)
%\eeaa
%for all $x\in [c, d]$ and $|\epsilon|\leq c/2.$ Notice that, for the same range of $x$ and $\epsilon$, the functions %$\frac{x}{x+\epsilon}$ and $\frac{k+x}{k+x+\epsilon}$ are increasing or decreasing in $x$ (depending on the sign of $\epsilon$)  at %the same time. Thus, $\frac{\Gamma(x+\epsilon)}{\Gamma(x)}$ is between $\frac{\Gamma(c+\epsilon)}{\Gamma(c)}$ and %$\frac{\Gamma(d+\epsilon)}{\Gamma(d)}.$ It follows that
%\beaa
%\sup_{c\leq x \leq d}\Big|\log \frac{\Gamma(x+\epsilon)}{\Gamma(x)}\Big| \leq \Big|\log \frac{\Gamma(c+\epsilon)}{\Gamma(c)}\Big| + %\Big|\log \frac{\Gamma(d+\epsilon)}{\Gamma(d)}\Big| \to 0
%\eeaa
%as $\epsilon\to 0$ since $\Gamma(x)$ is continuous for all $x>0.$\ \ \ \ \ \ \ $\blacksquare$

\noindent\textbf{Acknowledgement}. We thank Drs. Danning Li and Xingyun Zeng and Professors Xue Ding, Feng Luo, Albert Marden, Yongcheng Qi and Yong Zhang very much for their helps in discussing and checking the mathematical proofs in this paper. We also thank the referees and the associate editor for their valuable comments which improve the paper significantly.

\baselineskip 12pt
\def\ref{\par\noindent\hangindent 25pt}

\end{document}